\documentclass{amsart}

\usepackage{amssymb}
\usepackage{eepic}
\usepackage{color}

\allowdisplaybreaks

\numberwithin{equation}{section}

\newtheorem{theorem}{Theorem}[section]
\newtheorem{lemma}[theorem]{Lemma}

\newtheorem{remark}[theorem]{Remark}

\newtheorem{TheoA}{Theorem A}
\newtheorem{TheoB}{Theorem B1}
\newtheorem{TheoBB}{Theorem B2}

\newcommand{\Z}{\mathbb{Z}}
\newcommand{\R}{\mathbb{R}}
\newcommand{\C}{\mathbb{C}}

\newcommand{\summ}{\sum\nolimits}

\def\1{\mathbf{1}}
\def\H{\mathcal{H}}

\def\M{\mathcal{M}}
\def\A{\mathcal{A}}
\def\RR{\mathcal{R}}

\newcommand{\dem}{\noindent {\bf Proof. }}
\newcommand{\sketch}{\noindent {\bf Sketch of the proof. }}

\newcommand{\demB}{\noindent {\bf Proof of Theorem B1. }}
\newcommand{\demBB}{\noindent {\bf Proof of Theorem B2. }}

\newcommand{\fin}{\hspace*{\fill} $\square$ \vskip0.2cm}

\addtocounter{tocdepth}{-1}

\begin{document}

\title[Large BMO spaces vs interpolation]{Large BMO spaces vs interpolation}

\author[Conde-Alonso, Mei, Parcet]
{Jose M. Conde-Alonso, Tao Mei, Javier Parcet}

\maketitle

\null

\vskip-45pt

\null

\begin{abstract}
In this paper we introduce a class of BMO spaces which interpolate with $L_p$ and are sufficiently large to serve as endpoints for new singular integral operators. More precisely, let $(\Omega, \Sigma, \mu)$ be a $\sigma$-finite measure space. Consider two filtrations of $\Sigma$ by successive refinement of two atomic $\sigma$-algebras $\Sigma_\mathrm{a}, \Sigma_\mathrm{b}$ having trivial intersection. Construct the corresponding truncated martingale BMO spaces. Then, the intersection seminorm only leaves out constants and we provide a quite flexible condition on $(\Sigma_\mathrm{a}, \Sigma_\mathrm{b})$ so that the resulting space interpolates with $L_p$ in the expected way. In the presence of a metric $d$, we obtain endpoint estimates for Calder\'on-Zygmund operators on $(\Omega,\mu, d)$ under additional conditions on $(\Sigma_\mathrm{a}, \Sigma_\mathrm{b})$. These are weak forms of the \lq isoperimetric\rq${}$ and the \lq locally doubling\rq${}$ properties of Carbonaro/Mauceri/Meda which admit less concentration at the boundary. Examples of particular interest include densities of the form $e^{\pm |x|^\alpha}$ for any $\alpha > 0$ or $(1 + |x|^{\beta})^{-1}$ for any $\beta \gtrsim n^{3/2}$. A (limited) comparison with Tolsa's RBMO is also possible. On the other hand, a more intrinsic formulation yields a Calder\'on-Zygmund theory adapted to regular filtrations over $(\Sigma_\mathrm{a}, \Sigma_\mathrm{b})$ without using a metric. This generalizes well-known estimates for perfect dyadic and Haar shift operators. In contrast to previous approaches, ours extends to matrix-valued functions (via recent results from noncommutative martingale theory) for which only limited results are known and no satisfactory nondoubling theory exists so far. 
\end{abstract}

\addtolength{\parskip}{+1ex}

\null

\vskip-25pt

\null

\section*{{\bf Introduction}}

A BMO space is a set of functions which enjoy bounded mean oscillation in a certain sense. Both \lq\lq mean" and \lq\lq oscillation" can be measured in many different ways. Most frequently, we find BMO spaces referred to averages over balls in a metric measure space. In other notable scenarios, we may replace these averages by conditional expectations with respect to a martingale filtration or even by the action of a nicely behaved semigroup of operators. These more abstract formulations are known to be very useful in the lack of appropriate metrics. The relation between metric and martingale BMO spaces is well-understood for doubling spaces. In other words, when the measure of a ball in the given metric is comparable with the measure of its concentric dilations up to constants depending on the dilation factor but not on the chosen ball. Indeed, in this case the metric BMO is equivalent to a finite intersection of martingale BMO spaces constructed out of dyadic two-sided filtrations of atomic $\sigma$-algebras whose atoms look like balls, see \cite{C,GJ,HK,M1}. What is however more relevant is that any of these martingale BMO spaces satisfies the following fundamental properties:  
\begin{itemize}
\item[i)] Interpolation endpoint for the $L_p$ scale.

\item[ii)] John-Nirenberg inequalities and $\mathrm{H}_1 - \mathrm{BMO}$ duality.

\item[iii)] CZ extrapolation: $L_2$-boundedness $\Rightarrow$ $L_\infty \to \mathrm{BMO}$ boundedness.
\end{itemize}
Hence, these spaces yield at least as many endpoint estimates as the metric BMO. 

The main goal of this paper is to construct BMO spaces satisfying the properties stated above for a larger class of measures, and to explore the implications of it to provide new endpoint estimates. The first attempts in this direction \cite{MMNO,NTV1} culminated in the work of Tolsa \cite{To} on so-called RBMO spaces. These spaces enjoy the above mentioned properties for measures of polynomial growth. There are however a couple of open questions concerning Tolsa's construction. In first place Calder\'on-Zygmund extrapolation holds under a Lipschitz kernel condition instead of the more flexible H\"ormander condition. Second, only interpolation of operators has been studied but it seems to be unknown whether these spaces interpolate with the $L_p$ scale. These two problems were solved by Carbonaro, Mauceri and Meda \cite{CMM1,CMM2} for a different class of measures, based on similar results for the Gaussian measure on Euclidean spaces \cite{MM}. The properties they imposed lead to locally doubling measures with certain concentration behavior at the boundary. In both cases ---up to equivalence in the norm and additional conditions--- only doubling balls are used to measure the mean oscillation of the function.

We present an alternative approach to these questions. Martingale BMO spaces always satisfy conditions i) and ii) above, with independence of the existence of a metric in the underlying measure space. The third property however requires additional structure on our BMO spaces. Indeed, assume for a moment that we work with a two-sided filtration $(\Sigma_k)_{k \in \Z}$ of atomic $\sigma$-subalgebras of $\Sigma$ with corresponding conditional expectations $\mathsf{E}_{\Sigma_k}$. If $\mathbf{\Pi}$ denotes the union of atoms in our filtration, the corresponding martingale BMO norm is given by $$\|f\|_{\mathrm {BMO}} \, = \, \sup_{k \in \Z} \Big\| \mathsf{E}_{\Sigma_k} \big| f - \mathsf{E}_{\Sigma_{k-1}}f \big|^2 \Big\|_\infty^{\frac12}$$
which is larger than the following function-BMO norm $$\sup_{\mathrm{A} \in \mathbf{\Pi}} \Big( \frac{1}{\mu(\mathrm{A})} \int_{\mathrm{A}} \Big| f(w) - \frac{1}{\mu(\mathrm{A})} \int_{\mathrm{A}} f \, d\mu \Big|^2 d\mu(w) \Big)^{\frac12}.$$
Thus, if we admit from \cite{CMM1,To}Ê that extrapolation for (non-local) Calder\'on-Zygmund operators imposes our atoms to be doubling ---i.e. contained in a doubling ball of comparable measure or a union of at most $\mathrm{C}_0$ sets of this kind, see below--- we immediately find obstructions to construct filtrations satisfying this assumption for nondoubling spaces. We propose to consider a sort of intersection of two large BMO spaces as follows. Consider a $\sigma$-finite measure space $(\Omega, \Sigma, \mu)$ and two atomic $\sigma$-algebras $\Sigma_\mathrm{a}, \Sigma_\mathrm{b}$ of measurable sets in $\Sigma$ satisfying $\Sigma_\mathrm{a} \cap \Sigma_\mathrm{b} = \{\Omega, \emptyset\}$. Write $\mathrm{BMO}_j$ for any martingale BMO space over a filtration $(\Sigma_{jk})_{k \ge 1}$ with $\Sigma_{j1} = \Sigma_j$, then the semi-norm  
\begin{eqnarray*} 
\|f\|_{{\rm BMO}_{\Sigma_{\mathrm{ab}}}(\Omega)} = \max\Big\{ \big\| f - \mathsf{E}_{\Sigma_\mathrm{a}} f \big\|_{\mathrm{BMO}_\mathrm{a}} , \big\| f - \mathsf{E}_{\Sigma_\mathrm{b}} f \big\|_{\mathrm{BMO}_\mathrm{b}} \Big\} 
\end{eqnarray*} 
vanishes on constant functions precisely when $\Sigma_\mathrm{a} \cap \Sigma_\mathrm{b}$ is trivial. Let $$\mathrm{BMO}_{\Sigma_{\mathrm{ab}}}(\Omega) = \Big\{ f \in L_{\mathrm{loc}}^1(\Omega) \, \big| \ \|f  \|_{\mathrm{BMO}_{\Sigma_\mathrm{ab}}(\Omega)} < \infty \Big\} \big/ {\Bbb C}.$$ This settles a model of \lq large BMO spaces\rq${}$ which easily satisfy property ii) and leave some room for property iii). The problem reduces then to identify conditions on the pair $(\Sigma_\mathrm{a}, \Sigma_\mathrm{b})$ so that $\mathrm{BMO}_{\Sigma_{\mathrm{ab}}}(\Omega)$ interpolates with the $L_p$ scale. A standard argument shows that this is the case when $$\|f\|_{L_p^\circ(\Omega)} \, := \, \inf_{\mathrm{k} \in \C} \|f - \mathrm{k}\|_p \, \sim \, \max\Big\{ \|f - \mathsf{E}_{\Sigma_{\mathrm{a}}}f \|_p ,\|f - \mathsf{E}_{\Sigma_{\mathrm{b}}}f \|_p \Big\} \, =: \, \|f\|_{L^p_{\Sigma_{\mathrm{ab}}}(\Omega)}$$ for $2 \le p < \infty$, where 
\begin{eqnarray*}
L_p^\circ(\Omega) & = & L_p(\Omega, \Sigma, \mu)/\C, \\ [4pt] L^p_{\Sigma_{\mathrm{ab}}}(\Omega) & = & \big\{ f \in L_\mathrm{loc}^1(\Omega) \, \big| \ \|f\|_{L^p_{\Sigma_{\mathrm{ab}}}(\Omega)} < \infty \big\} / \C \\ & = & L_p(\Omega,\Sigma,\mu)/\Sigma_\mathrm{a} \bigwedge L_p(\Omega,\Sigma,\mu)/\Sigma_\mathrm{b}. 
\end{eqnarray*}
More precisely, we have an isomorphism $L_p^\circ(\Omega) \simeq L^p_{\Sigma_{\mathrm{ab}}}(\Omega)$. It should be mentioned that this isomorphism fails in general, even for the Lebesgue measure in $\R^n$ and many \lq natural\rq${}$ choices of pairs $(\Sigma_\mathrm{a}, \Sigma_\mathrm{b})$.Recall that $L_p^\circ(\Omega) = L_p(\Omega)$ for infinite measures. Note also that use use $\wedge$ and not $\cap$ since this space is not really an intersection; we shall also write $\mathrm{BMO}_{\Sigma_{\mathrm{ab}}}(\Omega) = \mathrm{BMO}_\mathrm{a}(\Omega)/\Sigma_\mathrm{a} \wedge \mathrm{BMO}_\mathrm{b}(\Omega)/\Sigma_\mathrm{b}$. To formulate a sufficient condition on $(\Sigma_\mathrm{a}, \Sigma_\mathrm{b})$ for $L_p^\circ(\Omega) \simeq L^p_{\Sigma_{\mathrm{ab}}}(\Omega)$, let $\Pi_j$ be the set of atoms in $\Sigma_j$. When $\mu(\Omega) < \infty$ we shall consider two distinguished atoms $(A_0, B_0) \in \Pi_\mathrm{a} \times \Pi_\mathrm{b}$, while for $\mu$ not finite we take $A_0 = B_0 = \emptyset$ for notation consistency. Given $(A,B) \in \Pi_\mathrm{a} \times \Pi_\mathrm{b}$, set $$R_{A} \, = \, \Big\{ B' \in \Pi_\mathrm{b} \, \big| \, \mu(A \cap B') > 0 \Big\} \quad \mbox{and} \quad R_{B} \, = \, \Big\{ A' \in \Pi_\mathrm{a} \, \big| \, \mu(A' \cap B) > 0 \Big\}.$$ We will write $|R_A|$ and $|R_B|$ for the cardinality of these sets. The following is the main result of this paper, where we establish a condition on $(\Sigma_\mathrm{a}, \Sigma_\mathrm{b})$ which suffices to make intersections and quotients commute in $L_p$ as described above. We will say that $(\Sigma_\mathrm{a},\Sigma_\mathrm{b})$ is an \emph{admissible covering} of $(\Omega, \Sigma,\mu)$ when $\Sigma_\mathrm{a} \cap \Sigma_\mathrm{b} = \{\Omega, \emptyset \}$ and $$\min \left\{ \sup_{A \in \Pi_\mathrm{a} \setminus \{A_0\}} \sum_{B \in R_A} |R_B| \, \frac{\mu(A \cap B)^2}{\mu(A) \mu(B)} \, , \, \sup_{B \in \Pi_\mathrm{b} \setminus \{B_0\}} \sum_{A \in R_B} |R_A| \, \frac{\mu(A \cap B)^2}{\mu(A) \mu(B)} \right\} \, < \, 1.$$ 

\begin{TheoA}
Let $(\Omega,\Sigma,\mu)$ be a $\sigma$-finite measure space equipped with an admissible covering $(\Sigma_\mathrm{a}, \Sigma_\mathrm{b})$. \hskip-4pt Then, for each $2 \le p < \infty$ there exists a constant $c_p$ depending only on $p$ and the admissible covering such that $$L_p^\circ(\Omega) \simeq L^p_{\Sigma_\mathrm{ab}}( \Omega)$$ with constant $c_p$. Moreover, we have the desired   complex interpolation result $$\big[ \mathrm{BMO}_{\Sigma_{\mathrm{ab}}}(\Omega), L_1^\circ (\Omega) \big]_{1/q} \simeq_{c_q} L_q^\circ(\Omega) \qquad (1 < q < \infty),$$ with $\mathrm{BMO}_{\Sigma_{\mathrm{ab}}}(\Omega)$ defined as above for any two martingale $\mathrm{BMO}$ spaces over  $(\Sigma_\mathrm{a}, \Sigma_\mathrm{b})$.
\end{TheoA}

The first assertion fails for $p = 1,\infty$ and it seems to be false for $p < 2$. On the other hand, both John-Nirenberg inequalities and $\mathrm{H}_1 - \mathrm{BMO}$ duality are easily formulated for these spaces. Therefore,  we shall focus in what follows on condition iii). Calder\'on-Zygmund extrapolation means that under certain mild smoothness condition on the kernel, $L_2$-boundedness yields $L_p$-boundedness for $1 < p < \infty$. As usual, we handle it by providing an endpoint estimate for interpolation. Let $d$ be a metric on $\Omega$ and denote by $\alpha \mathrm{B}$ the $\alpha$-dilation of a ball $\mathrm{B}$. We impose the standard H\"ormander kernel condition $$\sup_{\mathrm{B} \, d\mathrm{-ball}} \, \sup_{z_1, z_2 \in \mathrm{B}} \, \int_{\Omega \setminus \alpha \mathrm{B}}  \big| k(z_1,x) - k(z_2,x) \big| + \big| k(x,z_1) - k(x,z_2) \big| \, d\mu(x) \, < \, \infty.$$ Define a CZO on $(\Omega,\mu,d)$ as any linear map $T$ satisfying the following properties: 
\begin{itemize}
\item $T$ is well-defined and bounded on $L_2(\Omega)$.

\item The kernel representation for any $f \in \mathcal{C}_c(\Omega)$ $$Tf(x) \, = \, \int_\Omega k(x,y) f(y) \, d\mu(y) \quad \mbox{holds for} \quad x \notin \mathrm{supp} f$$ and some kernel $k: \Omega \times \Omega \setminus \Delta \to \C$ satisfying the H\"ormander condition.
\end{itemize}
Given $\mathrm{C}_0 > 0$, a $\Sigma$-measurable set $A$ will be called $(\mathrm{C}_0, \alpha, \beta)$-doubling when it is the union of at most $\mathrm{C}_0$ sets which are contained in $(\alpha,\beta)$-doubling balls of comparable measure up to the constant $\mathrm{C}_0$. 

\begin{TheoB}
Let $(\Sigma_\mathrm{a},\Sigma_\mathrm{b})$ be an admissible covering of $(\Omega, \Sigma, \mu)$. Assume that $(\Omega,\Sigma,\mu)$ admits regular filtrations $(\Sigma_{jk})_{k \ge 1}$ by successive refinement of $\Sigma_{j1} = \Sigma_j$ for $j =\mathrm{a},\mathrm{b}$ and that each atom in $\Sigma_{jk}$ is $(\mathrm{C}_0, \alpha, \beta)$-doubling for certain absolute constants $\mathrm{C}_0, \alpha, \beta > 0$. Construct the spaces $\mathrm{BMO}_{\Sigma_{\mathrm{ab}}}(\Omega)$ which are defined over these filtrations. Then, every Calder\'on-Zygmund operator extends to a bounded map $L_\infty(\Omega) \to \mathrm{BMO}_{\Sigma_{\mathrm{ab}}}(\Omega)$, and $L_p(\Omega) \to L_p(\Omega)$ for $1 < p < \infty$.
\end{TheoB}

\noindent A few illustrations of Theorem B1 are the following: 
\begin{itemize}
\item \textbf{Doubling case.} Theorem B1 recovers Calder\'on-Zygmund extrapolation on homogeneous spaces $(\Omega, \mu, d)$. We shall construct explicit pairs $(\Sigma_\mathrm{a}, \Sigma_\mathrm{b})$ and martingale filtrations satisfying our assumptions.

\vskip5pt

\item \textbf{Polynomial growth.} Given any $(\Omega, \mu,d)$ with polynomial growth, it is not difficult to construct atomic $\sigma$-algebras composed uniquely of doubling atoms, even giving admissible coverings. Under the existence of filtrations based on $(\Sigma_\mathrm{a}, \Sigma_\mathrm{b})$ and composed of doubling atoms ---regular or not--- we may prove that Tolsa's RBMO sits inside our $\mathrm{BMO}_{\Sigma_{\mathrm{ab}}}(\Omega)$. This condition seems unfortunately a restrictive limit in Theorem B1. However, it can be checked in some concrete scenarios like for $$d\mu(x) = \frac{dx}{1 + |x|^{\beta}} \quad \mbox{with} \quad \beta \gtrsim n^{\frac32}$$ in $\R^n$ equipped with the Euclidean metric. Note that $\mu$ is doubling for $\beta < n$. The key advantage with respect to Tolsa's approach is that we only need to impose H\"ormander kernel smoothness instead of stronger Lipschitz conditions. This was also achieved by Carbonaro-Mauceri-Meda for another family of measures (see below) but not for the measures considered above since they are drastically less concentrated at the boundary for any $\beta$. 

\vskip5pt

\item \textbf{Concentration at the boundary.} Carbonaro-Mauceri-Meda proved that when $(\Omega,\mu,d)$ is locally doubling and the measure concentrates at the boundary of open sets  in a certain sense ---together with a purely metric condition that does not play any role here--- a BMO space satisfying i), ii) and iii) is possible. Their main examples in $\R^n$ with a weighted Euclidean metric were $d\mu(x) = e^{\pm |x|^\alpha} dx$ and $\alpha > 1$. The exponentially decreasing ones behave in some sense like the Gaussian measure, which was studied a few years before by Mauceri-Meda. It is of polynomial growth, so that the kernel smoothness condition was the main advantage with respect to Tolsa's approach. The exponentially increasing ones are not of polynomial growth. In this paper we shall remove their condition $\alpha > 1$.  
\end{itemize}

In the literature, we find other families of operators ---with no need of a metric in the underlying space--- which are close to CZOs in spirit. Martingale transforms are the simplest ones, but local and much easier to bound. Non-local models include the so-called perfect dyadic CZOs and most notably Haar shift operators, which include prominent examples like the discrete Hilbert transform or dyadic paraproducts. In these cases, the H\"ormander kernel condition can be replaced by $$\sup_{\mathrm{Q} \, \mathrm{dyadic \, cube}} \, \sup_{z_1, z_2 \in \mathrm{Q}} \, \int_{\Omega \setminus \widehat{\mathrm{Q}}} \big| k(z_1,x) - k(z_2,x) \big| + \big| k(x,z_1) - k(x,z_2) \big| \, d\mu(x) \, < \, \infty,$$ where $\widehat{\mathrm{Q}}$ denotes the dyadic father of $\mathrm{Q}$. Our BMO spaces allow to further replace dyadic cubes in dyadically doubling measure spaces ---see \cite{LMP} for recent progress on more general measures in this direction---  by more general atoms. Namely, assume $(\Sigma_\mathrm{a}, \Sigma_\mathrm{b})$ gives an admissible covering of $(\Omega, \Sigma, \mu)$. Consider regular filtrations of atomic $\sigma$-algebras $(\Sigma_{jk})_{k \ge 1}$ with $\Sigma_{j1} = \Sigma_j$ for $j=\mathrm{a},\mathrm{b}$. Let us write $\Pi_{jk}$ for the family of atoms in the atomic $\sigma$-algebra $\Sigma_{jk}$ and set $\mathbf{\Pi}_j = \cup_{k \ge 1} \Pi_{jk}$. Then, consider the following H\"ormander-type kernel condition where the former role of the metric $d$ is replaced by the shape of our atoms in $\mathbf{\Pi} = \mathbf{\Pi}_\mathrm{a} \cup \mathbf{\Pi}_\mathrm{b}$ $$\sup_{A \in \mathbf{\Pi}} \, \sup_{z_1, z_2 \in A} \, \int_{\Omega \setminus \widehat{A}} \big| k(z_1,x) - k(z_2,x) \big| + \big| k(x,z_1) - k(x,z_2) \big| \, d\mu(x) \, < \, \infty.$$ Again, $\widehat{A}$ denotes the minimal atom in the filtration of $A$ which contains $A$ properly unless there is no such atom, in which case we pick $\widehat{A}=A$. If we replace the H\"ormander condition by this one, we obtain another class of \lq atomic\rq${}$ CZOs which will be denoted in what follows ACZO. 

\begin{TheoBB}
Let $(\Sigma_\mathrm{a},\Sigma_\mathrm{b})$ be an admissible covering of $(\Omega, \Sigma, \mu)$. Assume in addition that $(\Omega,\Sigma,\mu)$ admits regular filtrations $(\Sigma_{jk})_{k \ge 1}$ by successive refinement of $\Sigma_{j1} = \Sigma_j$ for $j=\mathrm{a},\mathrm{b}$. Construct the spaces $\mathrm{BMO}_{\Sigma_{\mathrm{ab}}}(\Omega)$ which are defined over these filtrations. Then, every $\mathrm{ACZO}$ extends to a bounded map from $L_\infty(\Omega)$ to $\mathrm{BMO}_{\Sigma_{\mathrm{ab}}}(\Omega)$. The $L_p$-boundedness follows as in Theorem \emph{B1}.
\end{TheoBB}

An advantage of Theorem B2 is that our kernel conditions are flexible since we may carefully choose $(\Sigma_a,\Sigma_b)$ and the regular filtrations according to the concrete singular integral operator. It is worth mentioning that every $\sigma$-finite (atomless if $\mu$ finite) measure space $(\Omega, \Sigma, \mu)$ has nontrivial admissible coverings. Of course the regularity of the filtration is a light form of \lq doublingness\rq${}$ needed to emulate the classical argument in this setting. We will also provide weaker estimates for pseudo-local operators when the filtrations are not regular.

In contrast to \cite{CMM1,CMM2,MM,To}, our approach extends to matrix-valued functions for which only limited results are known and no satisfactory nondoubling theory exists so far. In fact, this was our original motivation and the necessity of alternative arguments led to the results presented so far. We will postpone the discussion of the matrix-valued setting for the last section of this paper, which will allow those readers not familiar with noncommutative $L_p$ theory to isolate these results. 

Our results above give some insight on the relation between nondoubling and martingale BMO theories, see \cite{CP,JMP1} for other results along this line. In \cite{CP}, we adapt Tolsa's ideas to give an atomic block description of martingale $\mathrm{H}_1$. Semigroup BMO spaces are used in \cite{JMP1} to construct a Calder\'on-Zygmund theory which incorporates noncommutative measure spaces (von Neumann algebras) to the picture.


\noindent \textbf{Acknowledgement.} J.M. Conde-Alonso and J. Parcet are partially supported by the ERC StG-256997-CZOSQP, the Spanish grant MTM2010-16518 and by ICMAT Severo Ochoa Grant SEV-2011-0087 (Spain). T. Mei is partially supported by the NSF DMS-1266042.

\section{{\bf Admissible coverings and BMO spaces}} \label{Sect1}

In this section we recall some basic background around martingale BMO spaces and introduce our new class of BMO spaces. We will study standard properties of this class like the existence of admissible coverings, John-Nirenberg inequalities and $\mathrm{H}_1 - \mathrm{BMO}$ duality. The proof of Theorem A is more technical and will be postponed to Section \ref{SectMainTheorem}. 

\subsection{Martingale BMO spaces}

Let $(\Omega, \Sigma, \mu)$ be a $\sigma$-finite measure space and consider a filtration $(\Sigma_k)_{k \ge 1}$ of $\Sigma$. In other words, we have $\Sigma_k \subset \Sigma_{k+1}$ and the union of the spaces $L_\infty(\Omega, \Sigma_k, \mu)$ is weak-$*$ dense in $L_\infty(\Omega, \Sigma, \mu)$. Let $\mathsf{E}_{\Sigma_k}$ denote the conditional expectation onto $\Sigma_k$-measurable functions. Then, define the martingale $\mathrm{BMO}$  space associated to this filtration as the space of locally integrable functions $f: \Omega \to \C$ whose BMO norm below is finite $$\|f\|_{\mathrm{BMO}} \, = \, \sup_{k \ge 1} \Big\| \Big( \mathsf{E}_{\Sigma_k} \big| f - \mathsf{E}_{\Sigma_{k-1}}f \big|^2 \Big)^\frac12 \Big\|_\infty,$$ where we use the convention $\mathsf{E}_{\Sigma_0} f = 0$, see \cite{G}. Another expression for the norm is 
\begin{eqnarray*}
\|f\|_{\mathrm{BMO}} \!\!\! & = & \!\!\! \sup_{k \ge 1} \Big\| |df_k|^2 + \sum_{n > k} \mathsf{E}_{\Sigma_k} |df_n|^2 \Big\|_\infty^\frac12 \\ \!\!\! & \sim & \!\!\! \Big[ \sup_{k \ge 1} \Big\| \Big( \mathsf{E}_{\Sigma_k} \big| f - \mathsf{E}_{\Sigma_k} f \big|^2 \Big)^\frac12 \Big\|_\infty + \big\| \mathsf{E}_{\Sigma_1}f \big\|_\infty \Big] + \sup_{k > 1} \|df_k\|_\infty   
\end{eqnarray*}
where $df_k = \Delta_k f = \mathsf{E}_{\Sigma_k} f - \mathsf{E}_{\Sigma_{k-1}}f$. According to \cite{JJ}, $[\mathrm{BMO}, L_1(\Omega)]_{1/p} \simeq L_p(\Omega)$ for any filtration we pick. The bracket term in the right hand side above is called the martingale bmo norm of $f$ and it is closer to the standard expressions to measure the mean oscillation of a function. Namely, if $\Pi_k$ denotes the atoms in $\Sigma_k$ and $\mathbf{\Pi} = \cup_{k \ge 1} \Pi_k$, we deduce that 
\begin{eqnarray*}
\|f\|_{\mathrm{bmo}} \!\!\! & = & \!\!\! \sup_{k \ge 1} \Big\| \Big( \mathsf{E}_{\Sigma_k} \big| f - \mathsf{E}_{\Sigma_k} f \big|^2 \Big)^\frac12 \Big\|_\infty + \big\| \mathsf{E}_{\Sigma_1}f \big\|_\infty \\ \!\!\! & = & \!\!\! \sup_{A \in \mathbf{\Pi}} \Big( \frac{1}{\mu(A)} \int_A \Big| f(w) - \frac{1}{\mu(A)} \int_A f d\mu \Big|^2 d\mu(w) \Big)^\frac12 + \sup_{A \in \Pi_1} \Big| \frac{1}{\mu(A)} \int_A fd\mu \Big|.
\end{eqnarray*}
Of course, using a selected family of atoms makes $L_p$-interpolation fail in general for bmo. The extra term in BMO corrects this. This should be compared with the extra condition in the definition of Tolsa's RBMO. On the other hand, bmo spaces have good interpolation properties with little Hardy spaces $\mathrm{h}_p$. Namely, according to \cite{BCP} we have $[\mathrm{bmo}, \mathrm{h}_1]_{1/p} \simeq \mathrm{h}_p$ for any filtration and where $\mathrm{h}_p$ is the closure of the space of finite martingales in $L_p$ with respect to the norm $$\|f\|_{\mathrm{h}_p} \, = \, \Big\| \Big( \sum_{k \ge 1} \mathsf{E}_{\Sigma_{k-1}} |df_k|^2 \Big)^\frac12 \Big\|_p,$$ where this time the convention is $\mathsf{E}_{\Sigma_0} |df_1|^2 = |\mathsf{E}_{\Sigma_1} f|^2$. In contrast to other BMO spaces ---which always admit a null space--- the definitions above give norms, not seminorms. Paradoxically, we will need to quotient out certain spaces. Note that for $\Sigma_1$-measurable functions, the norms above coincide with the $L_\infty$ norm $$\|\mathsf{E}_{\Sigma_1}f\|_{\mathrm{BMO}} = \|\mathsf{E}_{\Sigma_1}f\|_{\mathrm{bmo}} = \|\mathsf{E}_{\Sigma_1}f\|_{L_\infty(\Omega)}.$$ If we define the seminorms
\begin{eqnarray*}
\|f\|_{\mathrm{bmo}}^\circ & = & \|f - \mathsf{E}_{\Sigma_1}f\|_{\mathrm{bmo}}, \\ [5pt] \|f\|_{\mathrm{BMO}}^\circ & = & \|f - \mathsf{E}_{\Sigma_1}f\|_{\mathrm{BMO}},
\end{eqnarray*}
we obtain complemented subspaces $\mathrm{BMO}_{\Sigma_1} = J_{\Sigma_1}(\mathrm{BMO})$ using the projection $J_{\Sigma_1} = id - \mathsf{E}_{\Sigma_1}$. Indeed, it is a simple exercise using Jensen's conditional inequality $|\mathsf{E}_{\Sigma_1} f |^2 \le \mathsf{E}_{\Sigma_1} |f|^2$, details are left to the reader. Since $J_{\Sigma_1}$ is also bounded on $\mathrm{h}_p$ and $L_p$, the previous interpolation results imply the following isomorphisms for $1 < p < \infty$ 
\begin{eqnarray*}
\big[ J_{\Sigma_1}(\mathrm{bmo}), J_{\Sigma_1}(\mathrm{h}_1(\Omega)) \big]_{1/p} & \simeq & J_{\Sigma_1}(\mathrm{h}_p(\Omega)), \\ \big[ J_{\Sigma_1}(\mathrm{BMO}), J_{\Sigma_1}(L_1(\Omega)) \big]_{1/p} & \simeq & J_{\Sigma_1}(L_p(\Omega)). 
\end{eqnarray*} 
Note that $J_{\Sigma_1}(L_p(\Omega)) \simeq L_p(\Omega,\Sigma,\mu)/\Sigma_1$ in the terminology of the Introduction.

\begin{remark} \label{LHS}
\emph{It is worth mentioning that Janson-Jones interpolation theorem \cite{JJ} holds for arbitrary filtrations. In particular, we could replace $(\Sigma_k)_{k \ge 1}$ by $(\Sigma_k)_{k \ge \mathrm{N}}$ for some large $\mathrm{N}$ and the latter BMO comes equipped with the norm $$\sup_{k \ge \mathrm{N}} \Big\| \Big( \mathsf{E}_{\Sigma_k} \big| f - \mathsf{E}_{\Sigma_k} f \big|^2 \Big)^\frac12 \Big\|_\infty + \big\| \mathsf{E}_{\Sigma_\mathrm{N}}f \big\|_\infty  + \sup_{k > \mathrm{N}} \|df_k\|_\infty.$$ When $\mathrm{N}$ is large enough, the middle term dominates the others and we get spaces which are closer and closer to $L_\infty(\Omega)$. On the contrary, when we quotient out the first $\sigma$-algebra by using the $J$-projections, it follows from the interpolation identities above that the starting $\sigma$-algebra affects significantly the interpolated space. This justifies in part our necessity to intersect two such spaces in this paper.}
\end{remark}

\subsection{BMO spaces for admissible coverings}

Let $(\Omega, \Sigma, \mu)$ be a $\sigma$-finite measure space and consider two atomic $\sigma$-algebras $\Sigma_\mathrm{a}, \Sigma_\mathrm{b}$ of measurable sets in $\Sigma$. Let $\Pi_j$ be the set of atoms in $\Sigma_j$ for $j = \mathrm{a}, \mathrm{b}$. When $\mu(\Omega) < \infty$, we shall consider two distinguished atoms $(A_0, B_0) \in \Pi_\mathrm{a} \times \Pi_\mathrm{b}$. If $\mu$ is not finite take $A_0 = B_0 = \emptyset$. Given $A \in \Pi_\mathrm{a}$, set $$R_{A} = \Big\{ B' \in \Pi_\mathrm{b} : \mu(A \cap B') > 0 \Big\}.$$ Define $R_{B}$ for $B \in \Pi_\mathrm{b}$ similarly. The pair $(\Sigma_\mathrm{a}, \Sigma_\mathrm{b})$ is called an \emph{admissible covering} of $(\Omega, \Sigma, \mu)$ when $\Sigma_\mathrm{a} \cap \Sigma_\mathrm{b} = \{\Omega, \emptyset\}$ and the inequality below holds: $$\min \left\{ \sup_{A \in \Pi_\mathrm{a} \setminus \{A_0\}} \sum_{B \in R_A} |R_B| \, \frac{\mu(A \cap B)^2}{\mu(A) \mu(B)} \, , \, \sup_{B \in \Pi_\mathrm{b} \setminus \{B_0\}} \sum_{A \in R_B} |R_A| \, \frac{\mu(A \cap B)^2}{\mu(A) \mu(B)} \right\} \, < \, 1.$$ Now, consider any pair of martingale filtrations $(\Sigma_{jk})_{k \ge 1}$ with $\Sigma_{j1} = \Sigma_j$ for $j = \mathrm{a}, \mathrm{b}$ and construct the corresponding martingale BMO spaces $\mathrm{BMO}_\mathrm{a}$ and $\mathrm{BMO}_\mathrm{b}$. As in the previous paragraph, we quotient out the $\Sigma_j$-measurable functions and set as we did in the Introduction
\begin{eqnarray*}
\mathrm{BMO}_{\Sigma_j}(\Omega) & = & J_{\Sigma_j}(\mathrm{BMO}_j), \\ \mathrm{BMO}_{\Sigma_\mathrm{ab}}(\Omega) & = & \mathrm{BMO}_{\Sigma_\mathrm{a}}(\Omega) \wedge \mathrm{BMO}_{\Sigma_\mathrm{b}}(\Omega)\\
&=&\big\{f\in L^1_{\mathrm{loc}}(\Omega) \, \big| \ \|f\|_{\mathrm{BMO}_{\Sigma_{\mathrm{ab}}}(\Omega)} < \infty \big\} / {\Bbb C}.
\end{eqnarray*}     

In the following we construct admissible coverings for $\sigma$-finite measure spaces. The procedure we employ is quite general. In concrete scenarios, other admissible coverings can be constructed enjoying additional properties as required in Theorems B1 and B2, these examples will be given later in this paper. 

\begin{remark}
\emph{The classical BMO on Euclidean spaces can be decomposed as an intersection of finitely many martingale BMO spaces, the amount of which depends on the dimension \cite{C,GJ,M1}. On the contrary, we just consider \lq\lq intersections" of two martingale BMOs. Note this makes our spaces larger and still amenable for interpolation, which gives some extra room to obtain endpoint estimates for singular integral operators. The main reason why this is possible is that our approach just relies on measure theoretic properties and does not rely on the geometry of the underlying space, as will become clear in the sequel.}
\end{remark}

\begin{lemma} \label{ACLemma}
Let $(\Omega, \Sigma, \mu)$ be a $\sigma$-finite measure space$\hskip1pt :$
\begin{itemize}
\item[i)] If $\mu(\Omega) = \infty$, it admits an admissible covering.

\item[ii)] If $\mu( \Omega) < \infty$ and $\mu$ is atomless, it admits an admissible covering.
\end{itemize}
\end{lemma}

\dem If $\mu(\Omega) = \infty$, pick $A_0 = \widetilde{A}_0 = B_0 = \widetilde{B}_0 = \emptyset$ and 
\begin{eqnarray*}
A_j & = & \widetilde{A}_j \setminus \widetilde{A}_{j-1}, \\ B_j & = & \widetilde{B}_j \setminus \widetilde{B}_{j-1},
\end{eqnarray*}
where $\emptyset \neq \widetilde{A}_1 \varsubsetneq \widetilde{B}_1 \varsubsetneq \widetilde{A}_2 \varsubsetneq \widetilde{B}_2 \varsubsetneq \widetilde{A}_3 \ldots$ are $\Sigma$-measurable sets chosen so that $$\min \left\{ \frac{\mu \big( \widetilde{B}_j \setminus \widetilde{B}_{j-1}\big)}{\mu \big( \widetilde{A}_j \big)}, \frac{\mu \big( \widetilde{A}_{j+1} \setminus \widetilde{A}_j \big)}{\mu \big( \widetilde{B}_j \big)} \right\} > \lambda > 4 \qquad \mbox{for all} \qquad j \ge 1.$$ It is at this point where we have used that $\mu(\Omega) = \infty$. Let $\Sigma_\mathrm{a}$ be the atomic $\sigma$-algebra generated by the atoms $(A_j)_{j \ge 1}$. Similarly, define $\Sigma_\mathrm{b} = \sigma \langle B_j: j \ge 1 \rangle$. It is clear by construction that $$\Sigma_\mathrm{a} \cap \Sigma_\mathrm{b} = \{\Omega, \emptyset\}.$$ On the other hand, $|R_B| = 2$ for every atom $B$ in $\Sigma_\mathrm{b}$. Therefore, it remains to show that $$\sup_{j \ge 1} \ \Big[ \frac{\mu(A_j \cap B_{j-1})^2}{\mu(A_j) \mu(B_{j-1})} + \frac{\mu(A_j \cap B_j)^2}{\mu(A_j) \mu(B_j)} \Big] \, < \, \frac12.$$ Note that the first summand above vanishes for $j=1$. The rest of terms are smaller than $1/\lambda$ according to our conditions, so that $\lambda > 4$ suffices. When $\mu(\Omega) < \infty$ we may assume that $\mu(\Omega)=1$ since renormalization does not affect our definition of admissible covering. We use again a \lq corona-type partition\rq${}$ $$\emptyset \neq \widetilde{A}_0 \varsubsetneq \widetilde{B}_0 \varsubsetneq \widetilde{A}_1 \varsubsetneq \widetilde{B}_1 \varsubsetneq \widetilde{A}_2 \ldots$$ satisfying $\mu(\widetilde{A}_0) = 1- \zeta$, $\mu(\widetilde{B}_0 \setminus \widetilde{A}_0) = \zeta (1-\zeta)$ and the relations below
\begin{eqnarray*}
\mu \big( \widetilde{A}_{j+1} \setminus \widetilde{B}_j \big) & = & \zeta \, \mu \big( \widetilde{B}_j \setminus \widetilde{A}_j \big), \\ \mu \big( \widetilde{B}_{j+1} \setminus \widetilde{A}_{j+1} \big) & = & \zeta \, \mu \big( \widetilde{A}_{j+1} \setminus \widetilde{B}_j \big),
\end{eqnarray*}
for $j \ge 0$. This is where we use the fact that $\mu$ has no atoms. Define $A_0 = \widetilde{A}_0$, $B_0 = \widetilde{B}_0$, $A_j = \widetilde{A}_j \setminus \widetilde{A}_{j-1}$ and $B_j = \widetilde{B}_j \setminus \widetilde{B}_{j-1}$ for $j \ge 1$. The $\sigma$-algebras $\Sigma_\mathrm{a}$ and $\Sigma_\mathrm{b}$ are the ones respectively generated by $(A_j)_{j \ge 0}$ and $(B_j)_{j \ge 0}$. In order to show that $\Omega = \cup_{j \ge 0} A_j = \cup_{j \ge 0} B_j$, let us prove that we have $$\sum_{j \ge 0} \mu(A_j) = \sum_{j \ge 0} \mu(B_j) = 1.$$ Indeed, if $j \ge 2$ we have 
\begin{eqnarray*}
\mu (A_j) \!\! & = & \!\! \mu \big( \widetilde{A}_j \setminus \widetilde{A}_{j-1} \big) \\ [3pt] \!\! & = & \!\! (1 + \zeta) \mu \big( \widetilde{B}_{j-1} \setminus \widetilde{A}_{j-1} \big) \ = \ \zeta (1 + \zeta) \mu \big( \widetilde{A}_{j-1} \setminus \widetilde{B}_{j-2} \big) \\ \!\! & = & \!\! \zeta (1 + \zeta) \Big[ \mu \big( \widetilde{A}_{j-1} \setminus \widetilde{A}_{j-2} \big) - \frac{1}{\zeta} \mu \big( \widetilde{A}_{j-1} \setminus \widetilde{B}_{j-2} \big) \Big] \ = \ \zeta^2 \mu(A_{j-1}).    
\end{eqnarray*}
Therefore, since $\mu(A_0) = 1 - \zeta$ and $\mu(A_1) = \zeta(1-\zeta^2)$, we deduce immediately that $\sum_{j \ge 0} \mu(A_j) = 1$. The sum $\sum_j \mu(B_j)$ also equals 1 since the two families are nested. The condition $\Sigma_\mathrm{a} \cap \Sigma_\mathrm{b} = \emptyset$ follows again by construction. Finally, since $|R_B|=2$ for all atoms $B=B_j$, it suffices one more time to prove that $$\sup_{j \ge 1} \ \Big[ \frac{\mu(A_j \cap B_{j-1})^2}{\mu(A_j) \mu(B_{j-1})} + \frac{\mu(A_j \cap B_j)^2}{\mu(A_j) \mu(B_j)} \Big] \, < \, \frac12.$$ According to our construction, the left hand side can be majorized by $$\frac{\mu(A_j \cap B_{j-1})^2}{\mu(A_j) \mu(B_{j-1})} + \frac{\mu(A_j \cap B_j)^2}{\mu(A_j) \mu(B_j)} \, \le \, \frac{\mu(\widetilde{B}_{j-1} \setminus \widetilde{A}_{j-1})}{\mu(B_{j-1})} + \frac{\mu(\widetilde{A}_j \setminus \widetilde{B}_{j-1})}{\mu(A_j)}.$$ On the other hand, arguing as before we may obtain the following identities $$\begin{array}{rclcrcl} \mu(A_j) & = & \zeta^{2(j-1)} (\zeta - \zeta^3), & & \mu \big( \widetilde{A}_j \setminus \widetilde{B}_{j-1} \big) & = & \zeta^{2j} (1 - \zeta), \\ \mu(B_{j-1}) & = & \zeta^{2(j-1)} (1 - \zeta^2), & & \mu(\widetilde{B}_{j-1} \setminus \widetilde{A}_{j-1}) & = & \zeta^{2(j-1)}(\zeta-\zeta^2). \end{array}$$ This gives a bound $2\zeta/(1+\zeta)$. It suffices for $\zeta < 1/3$. The proof is complete. \fin

\begin{remark}
\emph{All fully supported probability measures on $\R^n$ are nondoubling. In fact, this also holds for probability measures supported on unbounded sets. In particular, we hope Lemma \ref{ACLemma} together with Theorems B1 and B2 might open a door through further insight on Calder\'on-Zygmund theory for these measures. }
\end{remark}

\vskip-15pt 

\null

\subsection{John-Nirenberg inequalities, atomic $\mathrm{H}_1$ and duality}

In this paragraph we transfer some well-known properties of martingale BMO spaces to our new class of spaces. John-Nirenberg inequalities were formulated for the first time in \cite{JN}. Its analogue for martingale BMO spaces can be stated as follows: $$\sup_{k \ge 1} \sup_{A \in \Sigma_k} \, \frac{1}{\mu(A)} \, \mu \Big( A \cap \Big\{ \big| f - \mathsf{E}_{\Sigma_{k-1}}f \big| > \lambda \Big\} \Big) \, \lesssim \, \exp \Big( - \frac{c \lambda}{\|f\|_{\mathrm{BMO}}} \Big)$$ for all $\lambda > 0$, where the martingale $\mathrm{BMO}$ is constructed over the filtration $(\Sigma_k)_{k \ge 1}$ and we use the convention $\mathsf{E}_{\Sigma_0} f = 0$. The proof can be found in \cite{G}. An important consequence of this inequality is the $p$-invariance of the BMO norm. To be more precise, the martingale BMO norm admits the following equivalent expressions for any $0 < p < \infty$ $$\|f\|_{\mathrm{BMO}} \, \sim \, \sup_{k \ge 1} \Big\| \Big( \mathsf{E}_{\Sigma_k} \big| f - \mathsf{E}_{\Sigma_{k-1}} f \big|^p \Big)^{\frac{1}{p}} \Big\|_\infty.$$ If we replace $f$ by $J_{\Sigma_1}f = f - \mathsf{E}_{\Sigma_1}f$ in both inequalities, we immediately obtain the corresponding analogues for the BMO spaces which quotient out $\Sigma_1$-measurable functions, introduced above. Namely, the only difference is that we should read John-Nirenberg inequalities under the convention that $\mathsf{E}_{\Sigma_0} f = \mathsf{E}_{\Sigma_1}f$ and the $\mathrm{BMO}$ norm is given by $\| \ \|_{\mathrm{BMO}}^\circ$ instead. If we intersect two of these BMO spaces, we get John-Nirenberg type inequalities for our spaces $\mathrm{BMO}_{\Sigma_\mathrm{ab}}(\Omega)$ associated to an admissible covering $(\Sigma_\mathrm{a}, \Sigma_\mathrm{b})$ by taking again $\mathsf{E}_{\Sigma_0} f = \mathsf{E}_{\Sigma_1} f$
\begin{itemize}
\item $\displaystyle \|f\|_{\mathrm{BMO}_{\Sigma_\mathrm{ab}}(\Omega)} \, \sim \, \max_{j = \mathrm{a}, \mathrm{b}} \, \sup_{k \ge 1} \, \Big\| \Big( \mathsf{E}_{\Sigma_{jk}} \big| f - \mathsf{E}_{\Sigma_{j(k-1)}} f \big|^p \Big)^{\frac1p} \Big\|_\infty$, 

\item $\displaystyle \sup_{j = \mathrm{a}, \mathrm{b}} \sup_{\begin{subarray}{c} k \ge 1 \\ A \in \Sigma_{jk} \end{subarray}} \frac{1}{\mu(A)} \mu \Big( A \cap \Big\{ \big| f - \mathsf{E}_{\Sigma_{j(k-1)}} f \big| > \lambda \Big\} \Big) \, \lesssim \, \exp \Big( - \frac{c\lambda}{\|f\|_{\mathrm{BMO}_{\Sigma_\mathrm{ab}}(\Omega)}} \Big)$. 
\end{itemize}

Let us now consider $\mathrm{H}_1 - \mathrm{BMO}$ duality in our context. In the literature we find several equivalent descriptions of martingale $\mathrm{H}_1$ spaces via Doob's maximal function, martingale square function or conditional square function. Namely, $\mathrm{H}_1$ can be defined as the closure of the space of finite $L_1$ martingales with respect to any of the following norms $$\Big\| \sup_{k \ge 1} \big| \mathsf{E}_{\Sigma_k} f \big| \, \Big\|_1 \sim \Big\| \Big( \sum_{k \ge 1} |df_k|^2 \Big)^{\frac12} \Big\|_1 \sim \sum_{k \ge 1} \|df_k\|_1 + \Big\| \Big( \sum_{k \ge 1} \mathsf{E}_{\Sigma_{k-1}} |df_k|^2 \Big)^{\frac12} \Big\|_1.$$ We refer to Davis \cite{D} for the equivalences above and to Garsia \cite{G} for the duality theorem which claims that $\mathrm{H}_1^* \simeq \mathrm{BMO}$, a martingale analogue of Fefferman-Stein duality theorem. Let us now consider atomic descriptions of these spaces. The term \lq atom\rq${}$ unfortunately appears here in several settings ---$\sigma$-algebras, measures and Hardy spaces with different meanings--- but it will be clear which one is used from the context. Atomic descriptions are not possible for arbitrary $\mathrm{H}_1$ ---see \cite{CP} for an \lq atomic block\rq${}$ description both in the commutative/noncommutative settings---  but there are such results for $\mathrm{h}_1$ (defined above). A $\Sigma$-measurable function $a \in L_2(\Omega)$ is called an atom when there exists $k \ge 1$ and $A \in \Sigma_k$ with $$\mbox{supp} (a) \subset A, \qquad \mathsf{E}_{\Sigma_k} (a) = 0, \qquad \|a\|_2 \le \mu(A)^{-1/2}.$$ The atomic $\mathrm{h}_1$ is defined as the space of functions of the form $f = \summ_j \lambda_j a_j$ with $a_j$ atoms. The norm is the infimum of $\sum_j |\lambda_j|$ over all such possible expressions for the function $f$. This space is isomorphic to $\mathrm{h}_1$, see \cite{G}. In particular, it is also isomorphic to $\mathrm{H}_1$ when the filtration is regular. This will be enough for our purposes since we will only use $\mathrm{H}_1 - \mathrm{BMO}$ duality for regular filtrations. Now given two martingale filtrations $(\Sigma_{jk})_{k \ge 1}$ with $\Sigma_{j1} = \Sigma_j$ for $j = \mathrm{a}, \mathrm{b}$, let $\mathrm{H}_{1j}$ be the corresponding $\mathrm{H}_1$ spaces. Define $$\mathrm{H}_{\Sigma_{\mathrm{ab}}}^1(\Omega) = \Big\{ f \in L_1(\Omega) \ \big| \|f\|_{\rm H_1}=   \inf_{\begin{subarray}{c} f = f_1 + f_2 \\ \mathsf{E}_{\Sigma_\mathrm{a}}f_1=\mathsf{E}_{\Sigma_\mathrm{b}}f_2 = 0 \end{subarray}} \|f_1 \|_{\mathrm{H}_{1\mathrm{a}}} + \|f_2 \|_{\mathrm{H}_{1\mathrm{b}}}<\infty \Big\}.$$ Then, all the results above apply. In particular, we have $$\mathrm{H}_{\Sigma_{\mathrm{ab}}}^1(\Omega)^* \simeq \mathrm{BMO}_{\Sigma_{\mathrm{ab}}}(\Omega).$$

\section{{\bf Interpolation: Proof of Theorem A}}
\label{SectMainTheorem}

This section is entirely devoted to the proof of Theorem A. The argument is a bit lengthy, so that we have decided to divide it into several steps. We will assume that $\mu$ is a finite measure on $\Omega$ ---normalized so that $\mu(\Omega)=1$--- since this case is more technical. The slight modifications needed for the nonfinite case will be explained in the last step of the proof. 

\vskip3pt

\noindent \textbf{1. Intersection of quotients.} Let us first show that the interpolation result follows from the first assertion of Theorem A. Namely, given an admissible covering $(\Sigma_\mathrm{a}, \Sigma_\mathrm{b})$ of $(\Omega, \Sigma, \mu)$ and filtrations $(\Sigma_{jk})_{k \ge 1}$ with $\Sigma_{j1} = \Sigma_j$ for $j=\mathrm{a}, \mathrm{b}$, let $\mathrm{BMO}_j$ be the corresponding martingale BMO spaces. It is clear that 
\begin{eqnarray*}
\|f\|_{L_q^\circ(\Omega)} & = & \|f\|_{[ L_\infty^\circ(\Omega), L_1^\circ(\Omega) ]_{1/q}} \\ [2pt] & \gtrsim & \|f\|_{[ \mathrm{BMO}_{\Sigma_{\mathrm{ab}}}(\Omega), L_1^\circ(\Omega) ]_{1/q}} \\ [4pt] & \ge & \max_{j=\mathrm{a}, \mathrm{b}} \big\| f - \mathsf{E}_{\Sigma_j}f \big\|_{[ J_{\Sigma_j}(\mathrm{BMO}_j), J_{\Sigma_j}(L_1(\Omega))]_{1/q}} \\ &\simeq& \max_{j=\mathrm{a}, \mathrm{b}} \big\| f - \mathsf{E}_{\Sigma_j}f \big\|_{L_{\Sigma_{\mathrm{j}}}^q(\Omega)} \ = \ \|f\| _{L_{\Sigma_{\mathrm{ab}}}^q(\Omega). }
\end{eqnarray*}
For $q \ge 2$, this implies $$L_q^\circ (\Omega) \subset [\mathrm{BMO}_{\Sigma_{\mathrm{ab}}}(\Omega), L_1^\circ(\Omega) ]_{1/q} \subset L^q_{\Sigma_{\mathrm{ab}}}(\Omega).$$ Thus, the result follows from the isomorphism $L_q^\circ (\Omega) \simeq L^q_{\Sigma_{\mathrm{ab}}}(\Omega)$. The interpolation result for $1 < q <2$ follows from this and the well-known reiteration theorem \cite{BL}.

\vskip3pt

\noindent \textbf{2. Reduction to strict contractions.} The rest of the proof will be devoted to justify the first assertion of Theorem A. We claim that such isomorphism holds whenever we may find a constant $0 < c_p(\Sigma_{\mathrm{ab}}) < 1$ such that the following inequality holds for every mean-zero function $f \in L_p(\Omega)$ 
\begin{equation}
\label{strictcontractions}
\min \Big\{ \big\| \mathsf{E}_{\Sigma_\mathrm{a}} \mathsf{E}_{\Sigma_\mathrm{b}} f \big\|_p,  \big\| \mathsf{E}_{\Sigma_\mathrm{b}} \mathsf{E}_{\Sigma_\mathrm{a}} f \big\|_p \Big\} \, \le \, c_p(\Sigma_{\mathrm{ab}}) \|f\|_p.
\end{equation}
Indeed, if $\mathsf{E} \phi = \int_\Omega \phi \, d\mu$, we first observe that 
\begin{eqnarray*}
\|\phi\|_{L_p^\circ(\Omega)} & \sim & \big\| \phi - \mathsf{E} \phi \big\|_p \ \sim \ \inf_{\mathrm{k} \in \C} \big\| f-\mathrm{k} \big\|_p, \\ \|\phi\|_{J_{\Sigma_j}(L_p(\Omega))} & \sim & \big\| \phi - \mathsf{E}_{\Sigma_j} \phi \big\|_p \ \sim \ \inf_{\varphi \, \Sigma_j-\mathrm{measurable}} \big\| \phi - \varphi \big\|_{p}. 
\end{eqnarray*}
Therefore, our goal in what follows is to show that $$\big\| \phi - \mathsf{E} \phi \big\|_p \, \sim \, \big\| \phi - \mathsf{E}_{\Sigma_\mathrm{a}} \phi \big\|_p + \big\| \phi - \mathsf{E}_{\Sigma_\mathrm{b}} \phi \big\|_p$$ for every $\phi \in L_p(\Omega)$. The lower estimate is trivial. For the upper estimate, we shall use \eqref{strictcontractions}. Assume that the minimum above is attained at the first term (say) and let $f = \phi - \mathsf{E} \phi$ be a mean-zero function. We then find 
\begin{eqnarray*}
\big\| \mathsf{E}_{\Sigma_\mathrm{a}} \mathsf{E}_{\Sigma_\mathrm{b}} f \big\|_p & \le & c_p(\Sigma_{\mathrm{ab}}) \|f\|_p \\ & \le & c_p(\Sigma_{\mathrm{ab}}) \Big[ \big\| f - \mathsf{E}_{\Sigma_\mathrm{a}} f \big\|_p + \big\| \mathsf{E}_{\Sigma_\mathrm{a}} \big( f - \mathsf{E}_{\Sigma_\mathrm{b}} f \big) \big\|_p + \big\| \mathsf{E}_{\Sigma_\mathrm{a}} \mathsf{E}_{\Sigma_\mathrm{b}} f \big\|_p \Big],
\end{eqnarray*}
which implies 
$$\big\| \mathsf{E}_{\Sigma_\mathrm{a}} \mathsf{E}_{\Sigma_\mathrm{b}} f \big\|_p \, \le \, \frac{c_p(\Sigma_{\mathrm{ab}})}{1 - c_p(\Sigma_{\mathrm{ab}})} \Big[ \big\| \phi - \mathsf{E}_{\Sigma_\mathrm{a}} \phi \big\|_p + \big\| \phi - \mathsf{E}_{\Sigma_\mathrm{b}} \phi \big\|_p \Big].$$ This inequality is all we need, since the upper estimate follows from it 
\begin{eqnarray*}
\big\| \phi - \mathsf{E} \phi \big\|_p & \le & \big\| \mathsf{E}_{\Sigma_\mathrm{a}} \phi - \mathsf{E} \phi \big\|_p + \big\| \phi - \mathsf{E}_{\Sigma_\mathrm{a}} \phi \big\|_p \\ [5pt] & \le & \big\| \mathsf{E}_{\Sigma_\mathrm{a}} \mathsf{E}_{\Sigma_\mathrm{b}} f \big\|_p + \big\| \mathsf{E}_{\Sigma_\mathrm{a}} \big( \phi - \mathsf{E}_{\Sigma_{\mathrm{b}}} \phi \big) \big\|_p + \big\| \phi - \mathsf{E}_{\Sigma_\mathrm{a}} \phi \big\|_p \\ & \le & \frac{1}{1 - c_p(\Sigma_{\mathrm{ab}})} \Big[ \big\| \phi - \mathsf{E}_{\Sigma_\mathrm{a}} \phi \big\|_p + \big\| \phi - \mathsf{E}_{\Sigma_\mathrm{b}} \phi \big\|_p \Big].
\end{eqnarray*}

\vskip3pt

\noindent \textbf{3. The case $p = 2$.} Recall that we are assuming for the moment that $\mu(\Omega)=1$ and in that case we may consider two distinguished atoms $(A_0, B_0) \in \Sigma_\mathrm{a} \times \Sigma_\mathrm{b}$. In accordance with the previous point, it suffices to show that $$\min \Big\{ \big\| \mathsf{E}_{\Sigma_\mathrm{a}} \mathsf{E}_{\Sigma_\mathrm{b}} f \big\|_2, \big\| \mathsf{E}_{\Sigma_\mathrm{b}} \mathsf{E}_{\Sigma_\mathrm{a}} f \big\|_2 \Big\} \, \le \, c_2(\Sigma_{\mathrm{ab}}) \|f\|_2$$ for some $0 < c_2(\Sigma_\mathrm{ab}) < 1$ and every mean-zero $f \in L_2(\Omega)$. We claim that this estimate follows if the same inequality holds for $\Sigma_j$-measurable functions which vanish on the corresponding distinguished atom. More precisely, it suffices to prove that one of the following conditions holds
\begin{itemize}
\item $\big\| \mathsf{E}_{\Sigma_\mathrm{a}} \phi_\mathrm{b} \big\|_2 \, \le \, c_2(\Sigma_{\mathrm{ab}}) \|\phi_\mathrm{b}\|_2$ for $\phi_\mathrm{b}$ $\Sigma_\mathrm{b}$-measurable with $\phi_\mathrm{b}(B_0)=0$,

\item $\big\| \mathsf{E}_{\Sigma_\mathrm{b}} \phi_\mathrm{a} \big\|_2 \, \le \, c_2(\Sigma_{\mathrm{ab}}) \|\phi_\mathrm{a} \hskip0.5pt \|_2$ for $\phi_\mathrm{a}$ \hskip0.5pt $\Sigma_\mathrm{a}$\hskip0.5pt-measurable with $\phi_\mathrm{a} \hskip0.5pt (A_0)=0$.
\end{itemize}
Indeed, assume the first condition holds and let $\phi_\mathrm{b} \in L_2(\Omega, \Sigma_\mathrm{b}, \mu)$ be mean-zero. Then 
\begin{eqnarray*}
\|\phi_\mathrm{b}\|_2^2 & = & \big\| \phi_\mathrm{b} - \phi_\mathrm{b}(B_0) \big\|_2^2 - |\phi_\mathrm{b}(B_0)|^2, \\ \| \mathsf{E}_{\Sigma_\mathrm{a}} \phi_\mathrm{b} \|_2^2 & = & \big\| \mathsf{E}_{\Sigma_\mathrm{a}}(\phi_\mathrm{b} - \phi_\mathrm{b}(B_0)) \big\|_2^2 - |\phi_\mathrm{b}(B_0)|^2. 
\end{eqnarray*}
Substracting and using the first condition, we get $$\|\phi_\mathrm{b}\|_2^2 - \| \mathsf{E}_{\Sigma_\mathrm{a}} \phi_\mathrm{b} \|_2^2 \, \ge \, (1 - c_2(\Sigma_{\mathrm{ab}})) \big\| \phi_\mathrm{b} - \phi_\mathrm{b}(B_0) \big\|_2^2 \, \ge \, (1 -  c_2(\Sigma_{\mathrm{ab}})) \|\phi_\mathrm{b}\|_2^2.$$ Rearranging we get $\|\mathsf{E}_{\Sigma_\mathrm{a}} \phi_\mathrm{b} \|_2 \le c_2(\Sigma_{\mathrm{ab}}) \|\phi_\mathrm{b}\|_2$. Therefore, given any mean-zero $f \in L_2(\Omega)$, we may define $\phi_\mathrm{b} = \mathsf{E}_{\Sigma_\mathrm{b}} f$ and deduce that $\|\mathsf{E}_{\Sigma_\mathrm{a}} \mathsf{E}_{\Sigma_\mathrm{b}} f\|_2 \le c_2(\Sigma_{\mathrm{ab}}) \|f\|_2$ as desired. Alternatively, if we use the second condition above, the roles of $\Sigma_\mathrm{a}$ and $\Sigma_\mathrm{b}$ are switched and we obtain the other sufficient inequality which is implicit in the minimum above. Thus we have reduced the proof to justify one of the two conditions above. It is at this point where our definition of admissible pair comes into play. Namely, we know that $$\min \Big\{ \sup_{A \in \Pi_\mathrm{a} \setminus \{A_0\}} \sum_{B \in R_A} |R_B| \, \frac{\mu(A \cap B)^2}{\mu(A) \mu(B)}, \sup_{B \in \Pi_\mathrm{b} \setminus \{B_0\}} \sum_{A \in R_B} |R_A| \, \frac{\mu(A \cap B)^2}{\mu(A) \mu(B)} \Big\} = c(\Sigma_{\mathrm{ab}})$$ for some $0 < c(\Sigma_{\mathrm{ab}}) < 1$. Let us assume (say) that the minimum above is attained by the first term and let $\phi_\mathrm{a}$ be a $\Sigma_\mathrm{a}$-measurable function in $L_2(\Omega)$ that vanishes on $A_0$. Then, if we write $\phi_\mathrm{a} = \sum_{A \neq A_0} \alpha_A \chi_{A}$, we have the following estimate 
\begin{eqnarray*}
\|\mathsf{E}_{\Sigma_\mathrm{b}} \phi_\mathrm{a}\|_2^2 & = & \sum_{A,A' \neq A_0} \overline{\alpha}_A \alpha_{A'} \sum_{B \in R_{A} \cap R_{A'}} \frac{\mu(A \cap B)}{\mu(B)^\frac12} \frac{\mu(A' \cap B)}{\mu(B)^\frac12} \\ [5pt] & \le & \sum_{A,A' \neq A_0} \frac12 \sum_{B \in R_{A} \cap R_{A'}} \Big( |\alpha_A|^2 \frac{\mu(A \cap B)^2}{\mu(B)} + |\alpha_{A'}|^2 \frac{\mu(A' \cap B)^2}{\mu(B)} \Big) \\ [5pt] & = & \hskip5pt \sum_{A \neq A_0} |\alpha_A|^2 \mu(A) \sum_{\begin{subarray}{c} A' \neq A_0 \\ R_{A} \cap R_{A'} \neq \emptyset \end{subarray}} \sum_{B \in R_{A} \cap R_{A'}} \frac{\mu(A \cap B)^2}{\mu(A) \mu(B)} \\ & = & \hskip5pt \sum_{A \neq A_0} |\alpha_A|^2 \mu(A) \sum_{B \in R_{A}} |R_B| \frac{\mu(A \cap B)^2}{\mu(A) \mu(B)}  \ \le \ c(\Sigma_{\mathrm{ab}}) \sum_{A \neq A_0} |\alpha_A|^2 \mu(A).
\end{eqnarray*}
The right hand side equals $c(\Sigma_{\mathrm{ab}}) \|\phi_\mathrm{a}\|_2^2$, so we obtain the second condition. The first one follows when the minimum in our definition of admissible covering is attained by the second term. This proves that the first assertion of Theorem A holds for finite measures and $p=2$. The case $p>2$ requires some preliminaries. 

\vskip3pt

\noindent \textbf{4. A mass absorption principle.} Let us consider a particular ordering of the atoms in $\Sigma_\mathrm{a}$ and $\Sigma_{\mathrm{b}}$. According to our assumption $\Sigma_\mathrm{a} \cap \Sigma_\mathrm{b} = \{\Omega, \emptyset\}$, we may order $\Pi_\mathrm{a}$ so that $\Pi_{\mathrm a}=\{A_1, A_2, \ldots \}$ and for each $m \ge 0$ there exists $B \in \Pi_\mathrm{b}$ such that $\mu(A_{m+1} \cap B)$ and $\mu(\cup_{s \le m} A_s \cap B)$ are both strictly positive. Similarly, we may order $\Pi_b$ satisfying the symmetric condition. Define the atomic $\sigma$-algebras 
\begin{eqnarray*}
\Sigma_{\mathrm{a}}(m) & = & \sigma \Big\langle \bigcup_{s=0}^m A_s, \{A_s\}_{s \ge m+1}  \Big\rangle, \\ \Sigma_{\mathrm{b}}(m) & = & \sigma \Big\langle \bigcup_{s=0}^m B_s, \{B_s\}_{s \ge m+1}  \Big\rangle.  
\end{eqnarray*}
In this step we will prove that 
\begin{eqnarray}\label{reduc}
\|f \|_{L^p_{\Sigma_{\rm ab}}(\Omega)}  \simeq \big\| f - \mathsf{E}_{\Sigma_\mathrm{a}(m)}f \big\|_{L_p(\Omega)} + \big\| f - \mathsf{E}_{\Sigma_\mathrm{b}(m)}f \big\|_{L_p(\Omega)},
\end{eqnarray}
for any $m \ge 1$ and $2 < p < \infty$. The constants may depend on $m,p$ and the covering $(\Sigma_\mathrm{a}, \Sigma_\mathrm{b})$. Indeed, since the result is trivial for $m=0$, we will proceed by induction and assume that the result holds for $m-1$. Moreover, the upper estimate is straightforward and by symmetry it suffices to show that $$\big\| f - \mathsf{E}_{\Sigma_{\mathrm{a}}(m)} f \big\|_p \, \lesssim \, \big\| f - \mathsf{E}_{\Sigma_{\mathrm{a}}(m-1)} f \big\|_p + \big\| f - \mathsf{E}_{\Sigma_{\mathrm{b}}}f \big\|_p.$$ Taking $A_0(m) = \cup_{s \le m} A_s$, let $f = f\chi_{A_0(m)} + f \chi_{\Omega \setminus A_0(m)} = f_1 + f_2$. Since it is clear that $\mathsf{E}_{\Sigma_\mathrm{a}(m)} f_2 = \mathsf{E}_{\Sigma_\mathrm{a}(m-1)} f_2$, we may concentrate only on $f_1$. The left hand side for $f_1$ can be written as 
\begin{eqnarray*}
\big\| f_1 - \mathsf{E}_{\Sigma_{\mathrm{a}}(m)} f_1 \big\|_p & = & \big\| \chi_{A_0(m)} \big( f - \mathsf{E}_{\Sigma_{\mathrm{a}}(m)} f \big) \big\|_p \\ & = & \|f\|_{L_p^\circ(A_0(m))} \ \sim \ \sup_{\begin{subarray}{c}\|g\|_{L_{p'}(A_0(m))} \le 1 \\ g \, \mathrm{mean-zero} \end{subarray}} \Big| \int_{A_0(m)} fg \, d\mu \Big|.
\end{eqnarray*} 
Approximating the right hand side up to $\varepsilon >0$ by some mean-zero $g_0$ in the unit ball of $L_{p'}(A_0(m))$, let $B$ be an atom in $\Sigma_\mathrm{b}$ satisfying that $\mu(A_0(m-1) \cap B)$ and $\mu(A_m \cap B)$ are strictly positive. Recall that this can be done by the specific enumeration of atoms we picked. Then define
\begin{eqnarray*}
g_1 & = & \chi_{A_m}g_0 - \frac{\chi_{A_m \cap B}}{\mu(A_m \cap B)} \int_{A_m} g_0 \, d\mu, \\ g_2 & = & \chi_{A_0(m-1)}g_0 - \frac{\chi_{A_0(m-1) \cap B}}{\mu(A_0(m-1) \cap B)} \int_{A_0(m-1)} g_0 \, d\mu, \\ g_3 & = & \frac{\chi_{A_0(m-1) \cap B}}{\mu(A_0(m-1) \cap B)} \int_{A_0(m-1)} g_0 \, d\mu + \frac{\chi_{A_m \cap B}}{\mu(A_m \cap B)} \int_{A_m} g_0 \, d\mu. 
\end{eqnarray*}
Obviously, $g_0 = g_1 + g_2 + g_3$ and each $g_j$ is mean-zero. Moreover, we have  
\begin{eqnarray*}
\|g_1\|_{L_{p'}(A_0(m))} & \le & \| \chi_{A_m}g_0 \|_{L_{p'}(A_0(m))} + \Big\| \frac{\chi_{A_m \cap B}}{\mu(A_m \cap B)} \int_{A_m} g_0 \, d\mu \Big\|_{L_{p'}(A_0(m))} \\ & \le & \Big( 1 + \frac{\mu(A_m)^{\frac{1}{p}}}{\mu(A_m \cap B)^{\frac{1}{p}}} \Big) \| g_0 \|_{L_{p'}(A_0(m))} \ \lesssim \ \| g_0 \|_{L_{p'}(A_0(m))} \ \le \ 1.
\end{eqnarray*}
Similar computations apply to $g_2$ and $g_3$. In summary, we obtain the estimate below, where we write $f_Q$ to denote the average of $f$ over a given measurable set $Q$
\begin{eqnarray*}
\lefteqn{\hskip-10pt \big\| f_1 - \mathsf{E}_{\Sigma_{\mathrm{a}}(m)} f_1 \big\|_p} \\ & \sim & \Big| \int_{A_0(m)} f g_0 \, d\mu \Big| \\ & \le & \Big| \int_{A_m} \big( f - f_{A_m} \big) g_1 \, d\mu \Big| \\ & + & \Big| \int_{A_0(m-1)} \big( f - f_{A_0(m-1)} \big) g_2 \, d\mu \Big| + \Big| \int_{B} \big( f - f_B \big) g_3 \, d\mu \Big| \\ [6pt] & \lesssim & \big\| \chi_{A_m} \big( f - f_{A_m} \big) \big\|_p + \big\| \chi_{A_0(m-1)} \big( f - f_{A_0(m-1)} \big) \big\|_p + \big\| \chi_{B} \big( f - f_{B} \big) \big\|_p \\ [10pt] & \lesssim & \big\| f - \mathsf{E}_{\Sigma_{\mathrm{a}}(m-1)} f \big\|_p + \big\| f - \mathsf{E}_{\Sigma_{\mathrm{b}}}f \big\|_p.
\end{eqnarray*}
This completes the proof of the norm-equivalence \eqref{reduc}.

\vskip3pt

\noindent \textbf{5. The case $p > 2$.} We now complete the proof of Theorem A for probability measures. According to \eqref{reduc}, it suffices to show that there exists $0 < c_p(\Sigma_{\mathrm{ab}}) < 1$ and $m = m(p) \ge 1$ such that the following estimate holds for any mean-zero function $f \in L_p(\Omega)$ $$\min \Big\{ \big\| \mathsf{E}_{\Sigma_\mathrm{a}(m)} \mathsf{E}_{\Sigma_\mathrm{b}(m)} f \big\|_p, \big\| \mathsf{E}_{\Sigma_\mathrm{b}(m)} \mathsf{E}_{\Sigma_\mathrm{a}(m)} f \big\|_p \Big\} \, \le \, c_p(\Sigma_{\mathrm{ab}}) \|f\|_p.$$ Pick $m = m(p)$ as the smallest possible value of $m$ satisfying $$\min \Big\{ \mu(A_0(m)), \mu(B_0(m)) \Big\} \, > \, \max \Big\{ \Big( \frac{2 \cdot 4^{-p}}{1 - 2 \cdot 4^{-p}} \Big)^{\frac{1}{p-1}}, \big(1 - 4^{-p} \big)^{\frac{1}{p}} \Big\}$$ and $\varepsilon = \varepsilon (p) > 0$ small enough so that $$\big( 1 - 2 \cdot 4^{-p} \big)^\frac12 \, \le \, \big( 1 - 4^{-p} \big)^{\frac{1}{2p}} (1 - \varepsilon^3)^{\frac12} - \varepsilon^{\frac32}.$$ Since $L_p(\Omega) \subset L_2(\Omega)$, we know from Step 3 that $f$ always satisfies the above inequality for $p=2$. Assume that the minimum for $p=2$ is attained (say) at the first term so that $\| \mathsf{E}_{\Sigma_\mathrm{a}(m)} \mathsf{E}_{\Sigma_\mathrm{b}(m)} f \|_2 \le c_2(\Sigma_{\mathrm{ab}}) \|f\|_2$. Recall that $\mathsf{E}f = \int_\Omega f d\mu$. When $$\mathsf{E}(|\mathsf{E}_{\Sigma_{\mathrm{b}}(m)}f|^{\frac{p}{2}})^2 \, < \, (1 - \varepsilon^3) \|\mathsf{E}_{\Sigma_{\mathrm{b}}(m)}f\|_p^p$$ we proceed as follows 
\begin{eqnarray*}
\lefteqn{\hskip-15pt \big\| \mathsf{E}_{\Sigma_\mathrm{a}(m)} \mathsf{E}_{\Sigma_\mathrm{b}(m)} f \big\|_p^p} \\ [3pt] \!\! & \le & \!\! \big\| \mathsf{E}_{\Sigma_\mathrm{a}(m)} |\mathsf{E}_{\Sigma_\mathrm{b}(m)} f|^{\frac{p}{2}} \big\|_2^2 \hskip1pt - \hskip1pt \big( \mathsf{E} |\mathsf{E}_{\Sigma_{\mathrm{b}}(m)}f|^{\frac{p}{2}} \big)^2 + \big( \mathsf{E} |\mathsf{E}_{\Sigma_{\mathrm{b}}(m)}f|^{\frac{p}{2}} \big)^2 \\ [3pt] \!\! & = & \!\! \big\| \mathsf{E}_{\Sigma_\mathrm{a}(m)} \, \big( \, |\mathsf{E}_{\Sigma_{\mathrm{b}}(m)}f|^{\frac{p}{2}} - \mathsf{E} |\mathsf{E}_{\Sigma_{\mathrm{b}}(m)}f|^{\frac{p}{2}} \, \big) \big\|_2^2 \hskip1pt + \big( \mathsf{E} |\mathsf{E}_{\Sigma_{\mathrm{b}}(m)}f|^{\frac{p}{2}} \big)^2 \\ [3pt] \!\! & \le & \!\! \hskip3pt c_2^2(\Sigma_{\mathrm{ab}}) \, \big\| \, |\mathsf{E}_{\Sigma_{\mathrm{b}}(m)}f|^{\frac{p}{2}} \, - \, \mathsf{E} |\mathsf{E}_{\Sigma_{\mathrm{b}}(m)}f|^{\frac{p}{2}} \, \big\|_2^2 + \big( \mathsf{E} |\mathsf{E}_{\Sigma_{\mathrm{b}}(m)}f|^{\frac{p}{2}} \big)^2 \\ [5pt] \!\! & \le & \!\! \big[ c_2^2(\Sigma_\mathrm{ab}) + (1 - c_2^2(\Sigma_\mathrm{ab})) (1- \varepsilon^3) \big] \|f\|_p^p \, = \, c_p^p(\Sigma_\mathrm{ab}) \|f\|_p^p.
\end{eqnarray*}
If $\mathsf{E}(|\mathsf{E}_{\Sigma_{\mathrm{b}}(m)}f|^{\frac{p}{2}})^2 \ge (1 - \varepsilon^3) \|\mathsf{E}_{\Sigma_{\mathrm{b}}(m)}f\|_p^p$, then one can easily show that $$\big\| |\mathsf{E}_{\Sigma_{\mathrm{b}}(m)}f|^{\frac{p}{2}} - \mathsf{E} |\mathsf{E}_{\Sigma_{\mathrm{b}}(m)}f|^{\frac{p}{2}} \big\|_2^2 \, \le \, \varepsilon^3 \big\| \mathsf{E}_{\Sigma_{\mathrm{b}}(m)}f \big\|_p^p.$$ Now, decomposing $\mathsf{E}_{\Sigma_{\mathrm{b}}(m)}f = \mathsf{E}_{\Sigma_{\mathrm{b}}(m)}f(B_0(m)) \chi_{B_0(m)} + \mathsf{E}_{\Sigma_{\mathrm{b}}(m)}f \chi_{\Omega \setminus B_0(m)}$ we get
\begin{eqnarray*}
\lefteqn{\hskip-20pt \sqrt{\mu(B_0(m)) |\mathsf{E}_{\Sigma_{\mathrm{b}}(m)}f(B_0(m))|^p}} \\ & = & \big\| |\mathsf{E}_{\Sigma_{\mathrm{b}}(m)}f|^{\frac{p}{2}} \chi_{B_0(m)} \big\|_2 \\ [3pt] & \ge & \big\| \mathsf{E} |\mathsf{E}_{\Sigma_{\mathrm{b}}(m)}f|^{\frac{p}{2}} \chi_{B_0(m)} \big\|_2 - \big\| |\mathsf{E}_{\Sigma_{\mathrm{b}}(m)}f|^{\frac{p}{2}} - \mathsf{E} |\mathsf{E}_{\Sigma_{\mathrm{b}}(m)}f|^{\frac{p}{2}} \big\|_2 \\ [2pt] & \ge & \mu(B_0(m))^\frac12 \mathsf{E} | \mathsf{E}_{\Sigma_{\mathrm{b}}(m)}f|^{\frac{p}{2}} - \varepsilon^{\frac32} \hskip1pt \big\| \mathsf{E}_{\Sigma_{\mathrm{b}}(m)}f \big\|_p^{\frac{p}{2}} \\ & \ge & \Big[ \big( 1 - 4^{-p} \big)^{\frac{1}{2p}} (1 - \varepsilon^3)^{\frac12} - \hskip1pt \varepsilon^\frac32 \Big] \big\| \mathsf{E}_{\Sigma_{\mathrm{b}}(m)}f \big\|_p^{\frac{p}{2}} \\ & \ge & \big( 1 - 2 \cdot 4^{-p} \big)^\frac12 \big\| \mathsf{E}_{\Sigma_{\mathrm{b}}(m)}f \big\|_p^{\frac{p}{2}}.
\end{eqnarray*}
This also implies $$\big\| \mathsf{E}_{\Sigma_{\mathrm{b}}(m)}f \chi_{\Omega \setminus B_0(m)} \big\|_p^p \, \le \, 2 \cdot 4^{-p} \| \mathsf{E}_{\Sigma_{\mathrm{b}}(m)}f \|_p^p. $$ On the other hand, since $f$ is mean-zero we have $$\mathsf{E}_{\Sigma_{\mathrm{b}}(m)}f(B_0(m)) \mu(B_0(m)) + \mathsf{E} \big( \mathsf{E}_{\Sigma_{\mathrm{b}}(m)}f \chi_{\Omega \setminus B_0(m)} \big) = 0.$$ Rearranging and raising to the power $p$ then gives $$\mu(B_0(m))^p \big| \mathsf{E}_{\Sigma_{\mathrm{b}}(m)}f (B_0(m)) \big|^p \, \le \, \big\| \mathsf{E}_{\Sigma_{\mathrm{b}}(m)}f \chi_{\Omega \setminus B_0(m)} \big\|_p^p \, \le \, 2 \cdot 4^{-p} \big\| \mathsf{E}_{\Sigma_{\mathrm{b}}(m)}f \big\|_p^p.$$ Finally, combining our two estimates so far for $\mu(B_0(m))$ we obtain $$\mu(B_0(m)) \, \le \, \Big( \frac{2 \cdot 4^{-p}}{1 - 2 \cdot 4^{-p}} \Big)^{\frac{1}{p-1}}$$ which contradicts our choice of $m = m(p)$. This shows that $\mathsf{E}(|\mathsf{E}_{\Sigma_{\mathrm{b}}(m)}f|^{\frac{p}{2}})^2$ can not be larger than $(1 - \varepsilon^3) \|\mathsf{E}_{\Sigma_{\mathrm{b}}(m)}f\|_p^p$ and completes the proof in case the minimum for $p=2$ is attained at the first term. When the minimum is attained at the second term, a symmetric argument applies.  

\vskip3pt

\noindent \textbf{6. The nonfinite case.} When $\mu(\Omega) = \infty$ the proof of Theorem A is a bit simpler. In first place, note that $L_p^\circ(\Omega) = L_p(\Omega)$ in this case. In particular, the goal is to show that
\begin{eqnarray*}
L_q(\Omega) & \simeq & \big[ \mathrm{BMO}_{\Sigma_\mathrm{ab}}(\Omega), L_1(\Omega) \big]_{1/q}, \\ L_p(\Omega) & \simeq & L_p(\Omega, \Sigma, \mu)/\Sigma_\mathrm{a} \, \bigwedge \, L_p(\Omega, \Sigma, \mu)/ \Sigma_\mathrm{b}.
\end{eqnarray*}
Since $L_\infty(\Omega) \subset \mathrm{BMO}_{\Sigma_\mathrm{ab}}(\Omega)$, our argument in Step 1 can be easily adapted and interpolation follows from the second isomorphism above. To prove it, we follow essentially the same argument as for finite measures. Indeed, arguing as in Step 2 we see that it suffices to show that $$\min \Big\{ \big\| \mathsf{E}_{\Sigma_\mathrm{a}} \mathsf{E}_{\Sigma_\mathrm{b}} f \big\|_p,  \big\| \mathsf{E}_{\Sigma_\mathrm{b}} \mathsf{E}_{\Sigma_\mathrm{a}} f \big\|_p \Big\} \, \le \, c_p(\Sigma_{\mathrm{ab}}) \|f\|_p$$ for some constant $0 < c_p(\Sigma_{\mathrm{ab}}) < 1$ and every function $f \in L_p(\Omega)$. The only difference is that here it must hold for every $f$, not just mean-zero elements as in the finite case. The case $p=2$ is proved following Step 3. The fact that we do not assume $f$ to be mean-zero ---or ultimately to vanish at $A_0$ or $B_0$--- is compensated by our definition of admissible coverings, which does not consider distinguished atoms for infinite measures. Finally, once we know the case $p=2$ holds  ---for arbitrary functions, not only mean-zero ones--- we conclude that  $$\big\| \mathsf{E}_{\Sigma_\mathrm{a}} \mathsf{E}_{\Sigma_\mathrm{b}} f \big\|_p^p \, \le \, \big\| \mathsf{E}_{\Sigma_\mathrm{a}} \mathsf{E}_{\Sigma_\mathrm{b}} |f|^{\frac{p}{2}} \big\|_2^2 \, \le \, c_2^2(\Sigma_\mathrm{ab}) \big\| |f|^{\frac{p}{2}} \big\|_2^2 \, \le \, c_p^p(\Sigma_\mathrm{ab}) \|f\|_p^p$$ or a similar estimate for $\mathsf{E}_{\Sigma_\mathrm{b}} \mathsf{E}_{\Sigma_\mathrm{a}} f$. The proof of Theorem A is now complete. \fin 

\section{{\bf Calder\'on-Zygmund operators I}} \label{Sect3}

Let $(\Omega, \Sigma, \mu)$ be a measure space and consider a metric $d$ on $\Omega$. Assume that $\mu$ is $\sigma$-finite with respect to the metric topology. In this section we will be interested in Calder\'on-Zygmund operators on the metric measure space $(\Omega, \mu, d)$, as defined in the Introduction. More precisely, we prove Theorem B1 below and after that we shall illustrate this result with a few constructions of admissible coverings.  

\demB Our definition of CZO includes a symmetric H\"ormander kernel condition. This implies that the class of Calder\'on-Zygmund operators is closed under taking adjoints. In particular, the $L_p$-boundedness for $1 < p < 2$ can be deduced by duality from the case $p > 2$. On the other hand, according to Theorem A, the latter follows by interpolation if we can prove that any CZO extends to a bounded map $L_\infty(\Omega) \to \mathrm{BMO}_{\Sigma_\mathrm{ab}}(\Omega)$. Indeed, since $T$ is $L_2$-bounded Theorem A yields that $T: L_p(\Omega) \to L_p^\circ (\Omega).$ This is enough when the measure $\mu$ is infinite, since in that case $L_p(\Omega) = L_p^\circ (\Omega)$. When $\mu$ is finite we use $L_2$-boundedness once again together with H\"older's inequality to deduce that
\begin{eqnarray*}
\|Tf\|_p & \le & \|Tf - \mathsf{E}Tf\|_p + \mu(\Omega)^{\frac1p} |\mathsf{E}Tf| \\ & \lesssim & \|f\|_p + \mu(\Omega)^{\frac1p-\frac12} \|f\|_2 \ \lesssim \ \|f\|_p.   
\end{eqnarray*}
This completes the proof of our claim. Let us then prove the $L_\infty \to \mathrm{BMO}$ estimate. Consider an auxiliary BMO space which arises by averaging over the family of doubling balls in $(\Omega, \Sigma, \mu)$ $$\|f\|_{\mathrm{DBMO}} \, = \, \sup_{\begin{subarray}{c} \mathrm{B} \, d\mathrm{-ball} \\ \mathrm{doubling} \end{subarray}} \Big( \frac{1}{\mu(\mathrm{B})} \int_\mathrm{B} \Big| f(w) - \frac{1}{\mu(\mathrm{B})} \int_\mathrm{B} f \, d\mu \Big|^2 \, d\mu(w) \Big)^\frac12.$$ Following the standard argument, it is easily checked that $$T: L_\infty(\Omega) \to \mathrm{DBMO}.$$ Indeed, in first place we may observe as usual that we have the equivalence $$\|f\|_{\mathrm{DBMO}} \, \sim \, \sup_{\begin{subarray}{c} \mathrm{B} \, d\mathrm{-ball} \\ \mathrm{doubling} \end{subarray}} \inf_{\mathrm{k}_\mathrm{B} \in \C} \Big( \frac{1}{\mu(\mathrm{B})} \int_\mathrm{B} \big| f(w) - \mathrm{k}_\mathrm{B} \big|^2 \, d\mu(w) \Big)^\frac12.$$ Second, we decompose $f = f \chi_{\alpha \mathrm{B}} + f \chi_{\Omega \setminus \alpha \mathrm{B}} = \phi_{1\mathrm{B}} + \phi_{2\mathrm{B}}$ and pick the constant $\mathrm{k}_\mathrm{B}$ to be the average of $T\phi_{2\mathrm{B}}$ over $\mathrm{B}$. Then, we may estimate the norm of $Tf$ in $\mathrm{DBMO}$ by using the $L_2$-boundedness of $T$ for $T\phi_{1\mathrm{B}}$ and the H\"ormander kernel condition for $T\phi_{2\mathrm{B}}$. More precisely, we get
\begin{eqnarray*}
\|Tf\|_{\mathrm{DBMO}} & \le & \sup_{\begin{subarray}{c} \mathrm{B} \, d\mathrm{-ball} \\ \mathrm{doubling} \end{subarray}} \Big( \frac{1}{\mu(\mathrm{B})} \int_\mathrm{B} \big| T(f\chi_{\alpha \mathrm{B}})(w) \big|^2 \, d\mu(w) \Big)^\frac12 \\ & + & \Big( \frac{1}{\mu(\mathrm{B})} \int_\mathrm{B} \Big| T(f \chi_{\Omega \setminus \alpha \mathrm{B}})(w) - \frac{1}{\mu(\mathrm{B})} \int_\mathrm{B} T(f\chi_{\Omega \setminus \alpha \mathrm{B}}) \, d\mu \Big|^2 \, d\mu(w) \Big)^\frac12.
\end{eqnarray*}
Since we just use $(\alpha, \beta)$-doubling balls, the first term is dominated by $$\Big( \frac{\mu(\alpha \mathrm{B})}{\mu(\mathrm{B})} \Big)^{\frac12} \big\| T \big\|_{2 \to 2} \, \|f\|_\infty \, \lesssim \, \|f\|_\infty.$$ On the other hand, using the kernel representation of $T$ we may write 
\begin{eqnarray*}
\lefteqn{\hskip-25pt T(f \chi_{\Omega \setminus \alpha \mathrm{B}})(w) - \frac{1}{\mu(\mathrm{B})} \int_\mathrm{B} T(f\chi_{\Omega \setminus \alpha \mathrm{B}}) \, d\mu} \\ \hskip20pt & = & \frac{1}{\mu(\mathrm{B})} \int_\mathrm{B} \int_{\Omega \setminus \alpha \mathrm{B}} \big( k(w,\zeta) - k(\xi,\zeta) \big) f(\zeta) \, d\mu(\zeta) \, d\mu(\xi)
\end{eqnarray*}
for $w \in \mathrm{B}$. In particular, the last term above can be majorized by $\|f\|_\infty$ using H\"ormander condition for $k$. This proves the $L_\infty(\Omega) \to \mathrm{DBMO}$ boundedness of our CZO. Therefore, it suffices to show that $\mathrm{DBMO} \subset \mathrm{BMO}_{\Sigma_\mathrm{ab}}(\Omega)$. This follows from the following chain of inclusions $$\mathrm{DBMO} \, \subset \, \mathrm{bmo}_{\Sigma_\mathrm{a}} \wedge \mathrm{bmo}_{\Sigma_\mathrm{b}} \, \subset \, \mathrm{BMO}_{\Sigma_\mathrm{a}} \wedge \mathrm{BMO}_{\Sigma_\mathrm{b}} \, = \, \mathrm{BMO}_{\Sigma_\mathrm{ab}}(\Omega).$$ Let us recall in passing the terminology we are using, namely $$\mathrm{bmo}_{\Sigma_j} \, = \, J_{\Sigma_j}(\mathrm{bmo}_j) \qquad \mbox{and} \qquad \mathrm{BMO}_{\Sigma_j} \, = \, J_{\Sigma_j} (\mathrm{BMO}_j)$$ for $j = \mathrm{a}, \mathrm{b}$. Here $\mathrm{bmo}_j$ and $\mathrm{BMO}_j$ are the martingale bmo and BMO spaces constructed over the filtrations $(\Sigma_{jk})_{k \ge 1}$ described in the statement of Theorem B1. If $\mathbf{\Pi}_j$ denotes the atoms in such filtration, the norm in $\mathrm{bmo}_{\Sigma_j}$ is given by 
\begin{eqnarray*} 
\|f\|_{\mathrm{bmo}_{\Sigma_j}} & = & \sup_{k \ge 1} \Big\| \mathsf{E}_{\Sigma_{jk}} \big| f - \mathsf{E}_{\Sigma_{jk}} f \big|^2 \Big\|_\infty^\frac12 \\ & = & \sup_{A \in \mathbf{\Pi}_j} \Big( \frac{1}{\mu(A)} \int_A \Big| f(w) - \frac{1}{\mu(A)} \int_A f \, d \mu \Big|^2 d\mu(w) \Big)^\frac12.
\end{eqnarray*}
Now, since we assume that all atoms in $\mathbf{\Pi} = \mathbf{\Pi}_\mathrm{a} \cup \mathbf{\Pi}_\mathrm{b}$ are doubling, the seminorm above is majorized (up to absolute constants) by the seminorm in DBMO. As this holds for both $j = \mathrm{a}, \mathrm{b}$, we have proved the first inclusion. Now, for the second inclusion, we recall the seminorm in $\mathrm{BMO}_{\Sigma_j}$ $$\|f\|_{\mathrm{BMO}_{\Sigma_j}} \, = \, \sup_{k \ge 1} \Big\| \mathsf{E}_{\Sigma_{jk}} \big| f - \mathsf{E}_{\Sigma_{jk-1}}f \big|^2 \Big\|_\infty^\frac12$$ where $\mathsf{E}_{\Sigma_{j0}}f = \mathsf{E}_{\Sigma_{j1}}f$ since we quotient out $\Sigma_{j1}$-measurable functions. Note also that we are imposing the filtrations $(\Sigma_{jk})_{k \ge 1}$ to be regular. In other words, there exist absolute constants $c_j > 0$ such that $\mathsf{E}_{\Sigma_{jk}} |f| \le c_j \mathsf{E}_{\Sigma_{jk-1}}|f|$ for $j=\mathrm{a}, \mathrm{b}$ and $k \ge 1$. This yields the inequality $$\|f\|_{\mathrm{BMO}_{\Sigma_j}} \, \le \, c_j \|f\|_{\mathrm{bmo}_{\Sigma_j}}.$$ Thus, $\mathrm{BMO}_{\Sigma_{\mathrm{ab}}}(\Omega) \simeq \mathrm{bmo}_{\Sigma_\mathrm{a}} \wedge \mathrm{bmo}_{\Sigma_\mathrm{b}}$ for regular filtrations and we are done. \fin 

\begin{remark} \label{H1Rem}
\emph{Under the same assumptions, every $\mathrm{CZO}$ extends to a bounded map $$\mathrm{H}_{\Sigma_{\mathrm{ab}}}^1(\Omega) \to L_1(\Omega).$$ Indeed, it follows at once by duality and Theorem B1. Alternatively, since we need to work with regular filtrations, we may use the atomic description of $\mathrm{H}_{\Sigma_{\mathrm{ab}}}^1(\Omega)$ given in Section \ref{Sect1} from which an easy argument arises, details are left to the reader.}  
\end{remark}

\noindent In the following paragraphs we shall illustrate Theorem B1 with a few examples.

\subsection{Doubling case} \label{DCSubsection}

Admissible coverings fulfilling the assumptions in Theorem B1 can always be constructed on every doubling space, so that Calder\'on-Zygmund extrapolation for homogeneous spaces appears as a particular application of our approach. For clarity of the exposition, we shall just indicate how to construct such admissible coverings in $\R^2$ with the Lebesgue measure $m$ and the Euclidean metric, although a similar construction works in the general case. Let us pick $\mathrm{Q}_0 = [ -\frac12, \frac12 ] \times [ -\frac12, \frac12 ]$ the unit cube and set $\mathrm{Q}_s = 3^s \mathrm{Q}_0$ for $s \ge 1$. Consider the $\sigma$-algebras  $$\Sigma_\mathrm{a} \, = \, \sigma \big\langle A_s : s \ge 1 \big\rangle \quad \mbox{and} \quad \Sigma_\mathrm{b} \, = \, \sigma \big\langle B_s : s \ge 1 \big\rangle,$$ where $(A_1, B_1) \hskip1pt = \hskip1pt (\mathrm{Q}_0, \mathrm{Q}_1)$ and $(A_s,B_s) \hskip1pt = \hskip1pt (\mathrm{Q}_{2s-2} \setminus \mathrm{Q}_{2s-4} \hskip1pt , \hskip1pt \mathrm{Q}_{2s-1} \setminus \mathrm{Q}_{2s-3})$ for $s \ge 2$. Then it follows from the proof of Lemma \ref{ACLemma} that $(\Sigma_\mathrm{a}, \Sigma_\mathrm{b})$ is an admissible covering of the Euclidean space $(\R^2, m)$. 

\begin{picture}(360,200)(-170,-100)
\linethickness{1.5pt}

    \put(-9.1,-8.5){\line(1,0){18}}
    \put(-8.8,-7.8){\line(1,0){18}}
    \put(-8.8,-6.6){\line(1,0){18}}
    \put(-8.8,-5.4){\line(1,0){18}}
    \put(-8.8,-4.2){\line(1,0){18}}
    \put(-8.8,-3.0){\line(1,0){18}}
    \put(-8.8,-1.8){\line(1,0){18}}
    \put(-8.8, 0.6){\line(1,0){18}}
    \put(-8.8,4.2){\line(1,0){18}}
    \put(-8.8,3.0){\line(1,0){18}}
    \put(-8.8,1.8){\line(1,0){18}}
    \put(-8.8,5.4){\line(1,0){18}}
    \put(-8.8,7.8){\line(1,0){18}}
    \put(-8.8,6.6){\line(1,0){18}}
    \put(-8.8,0.0){\line(1,0){18}}
    \put(-8.8,1.0){\line(1,0){18}}

    \put(-8.5,-8.5){\line(0,1){17.7}}
    \put(-7.8,-8.5){\line(0,1){17.7}}
    \put(-6.6,-8.5){\line(0,1){17.7}}
    \put(-5.4,-8.5){\line(0,1){17.7}}
    \put(-4.2,-8.5){\line(0,1){17.7}}
    \put(-3.0,-8.5){\line(0,1){17.7}}
    \put(-1.8,-8.5){\line(0,1){17.7}}
    \put(0.6,-8.5){\line(0,1){17.7}}
    \put(-0.6,-8.5){\line(0,1){17.7}}
    \put(4.2,-8.5){\line(0,1){17.7}}
    \put(3.0,-8.5){\line(0,1){17.7}}
    \put(1.8,-8.5){\line(0,1){17.7}}
    \put(5.4,-8.5){\line(0,1){17.7}}
    \put(8.5,-9){\line(0,1){18.1}}
    \put(7.8,-8.5){\line(0,1){17.7}}
    \put(6.6,-8.5){\line(0,1){17.7}}
    \put(0.0,-8.5){\line(0,1){17.7}}
    \put(1.0,-8.5){\line(0,1){17.7}}
  
    
    \put(-82.5,-82.5){\line(0,1){165}}
    \put(-82.5,-82.5){\line(1,0){165}}
    \put(82.5,82.5){\line(0,-1){165}}
    \put(82.5,82.5){\line(-1,0){165}}

\linethickness{0.8pt}

    \put(-27.5,-27.5){\line(0,1){55}}
    \put(-27.5,-27.5){\line(1,0){55}}
    \put(27.5,27.5){\line(0,-1){55}}
    \put(27.5,27.5){\line(-1,0){55}}

\linethickness{.3pt}

    \put(82.5,9.16){\line(-1,0){165}}
    \put(82.5,27.5){\line(-1,0){165}}
    \put(82.5,45.83){\line(-1,0){165}}
    \put(82.5,64.16){\line(-1,0){165}}
    \put(82.5,-9.17){\line(-1,0){165}}
    \put(82.5,-27.5){\line(-1,0){165}}
    \put(82.5,-45.83){\line(-1,0){165}}
    \put(56,-64.16){\line(-1,0){138}}
    
    \put(-95,-82.5){\line(1,0){190}}
    \put(-82.5,-95){\line(0,1){190}}
    \put(82.5,95){\line(0,-1){190}}
    \put(95,82.5){\line(-1,0){190}}
    
    \put(9.16,82.5){\line(0,-1){165}}
    \put(27.5,82.5){\line(0,-1){165}}
    \put(45.83,82.5){\line(0,-1){165}}
    \put(64.16,82.5){\line(0,-1){138}}
    \put(-9.17,82.5){\line(0,-1){165}}
    \put(-27.5,82.5){\line(0,-1){165}}
    \put(-45.83,82.5){\line(0,-1){165}}
    \put(-64.16,82.5){\line(0,-1){165}}

\linethickness{.1pt}

    \put(60,-69){{\Large {\bf \emph{B}}$_{\hskip-2pt \mbox{{\small{\bf \emph{s}}}}}$}}
    \put(12,-20){{\small$A_s$}}

    \multiput(0,-82.5)(0,10){6}{\line(0,1){5}}
    \multiput(0,82.5)(0,-10){12}{\line(0,-1){5}}
    \multiput(-82.5,0)(10,0){12}{\line(1,0){5}}
    \multiput(82.5,0)(-10,0){6}{\line(-1,0){5}}
\end{picture}

\null

\vskip-35pt 

\null

\begin{center} \textsc{Figure I} \\ {\small The admissible covering and the second generation of one of the filtrations}
\end{center}

Next we define the filtrations $(\Sigma_{jk})_{k \ge 1}$ with $\Sigma_{j1} = \Sigma_j$ for $j = \mathrm{a}, \mathrm{b}$. Except for $A_1$ and $B_1$ ---which are ordinary cubes--- the atoms $A_s$ and $B_s$ ($s \ge 2$) are punctured cubes in which we remove a concentric cube with side-length $1/9$ times the side-length of the larger one. To define $\Sigma_{j2}$ for $j = \mathrm{a}, \mathrm{b}$, we break each $A_s, B_s$ into a disjoint union of $80$ equal cubes of side-length $1/9$ times the side-length of the original punctured cube ---i.e. all except for the one in the centre--- unless $s=1$ in which case we also pick the centre and get $81$ subcubes. The next generations are simpler. Indeed, since all our atoms in $\Sigma_{j2}$ are already cubes, we perform dyadic partitions in each of them to provide the next generations of our filtration. This procedure completely defines two filtrations respectively based on $\Sigma_\mathrm{a}$ and $\Sigma_\mathrm{b}$. It remains to check that these filtrations are regular and the atoms are doubling. The regularity constant is dominated uniformly by $81$ when $(k-1,k)=(1,2)$ and by $4$ otherwise. On the other hand, our atoms for $k=1$ are punctured cubes which are comparable to the corresponding unpunctured ones, which in turn are doubling with constant $4$. This proves that all conditions in Theorem B1 are satisfied. In the general case we just need to use Christ dyadic cubes and adapt our choice according to the finiteness or non-finiteness of $\mu$ as we did in Lemma \ref{ACLemma}.

\subsection{Polynomial growth}

Given $(\Omega, \mu, d)$ of polynomial growth with $\mu(\Omega) = \infty$ we may easily construct an admissible covering of $(\Omega, \mu)$ composed of doubling atoms. Indeed, the construction above can be easily modified using the existence of arbitrarily large doubling cubes centered at almost every point in the support of $\mu$, see \cite{To} for details. The main difficulty relies in the construction of martingale filtrations $(\Sigma_{jk})_{k \ge 1}$ satisfying the assumptions in Theorem B1. Note that, whenever that holds, we find $$\mathrm{RBMO} \, \subset \, \mathrm{DBMO} \, \subset \, \mathrm{BMO}_{\Sigma_{\mathrm{ab}}}(\Omega).$$ In particular, we deduce that $[\mathrm{RBMO}, L_1(\Omega)]_{1/q} \simeq L_q(\Omega)$ when this happens. As far as we know, such interpolation identities are new since Tolsa studied in \cite{To} interpolation of operators. Unfortunately, the construction of such filtrations seems to be a difficult task in the general case. For instance, the corona-type construction described above finds some obstructions when the measure $\mu$ is supported in Cantor like sets. Nevertheless, we may construct these filtrations in some other cases. Let us consider the following family of measures on $\R^n$ equipped with the Euclidean distance $$d\mu_\beta(x) \, = \, \frac{dx}{1 + |x|^\beta}.$$ 

\noindent These measures are nondoubling only for $\beta > n$. We will construct an admissible covering for $\beta \gtrsim n^{3/2}$ satisfying the hypotheses of Theorem B1 when $d$ is the Euclidean metric in $\R^n$. We will work with the equivalent measure $$d\nu_\beta(x) = \min \{1, |x|^{-\beta}\} \, dx$$ for convenience. Note that this does not affect the conclusions in Theorem B1.  

Pick $\mathrm{Q}_0 = [ - \lambda, \lambda ]^n$ with $\lambda> 1$ to be fixed and set $\mathrm{Q}_s = 2^s \mathrm{Q}_0$. Consider the $\sigma$-algebras $\Sigma_\mathrm{a} = \sigma \big\langle A_s : s \ge 1 \big\rangle$ and $\Sigma_\mathrm{b} = \sigma \big\langle B_s : s \ge 1 \big\rangle$ where $(A_0, B_0) \hskip1pt = \hskip1pt (\mathrm{Q}_0, \mathrm{Q}_1)$ and $(A_s,B_s) \hskip1pt = \hskip1pt (\mathrm{Q}_{2s} \setminus \mathrm{Q}_{2s-2} \hskip1pt , \hskip1pt \mathrm{Q}_{2s+1} \setminus \mathrm{Q}_{2s-1})$ for any $s \ge 1$. We clearly have $\Sigma_\mathrm{a} \cap \Sigma_\mathrm{b} = \{\R^n, \emptyset\}$ and $\max \{|R_A|, |R_B| \} \le 2$ for $(A,B) \in \Pi_\mathrm{a} \times \Pi_\mathrm{b}$ by construction of $(\Sigma_\mathrm{a}, \Sigma_\mathrm{b})$. Thus, it suffices to show that $$\sup_{\begin{subarray}{c} (A,B) \in \Pi_\mathrm{a} \times \Pi_\mathrm{b} \\ A \neq A_0 \end{subarray}} \frac{\nu_\beta(A \cap B)^2}{\nu_\beta(A)\nu_\beta(B)} \, < \, \frac12.$$ We will prove in fact the apparently stronger inequalities $$\sup_{s \ge 1} \frac{\nu_\beta(A_s \cap B_{s-1})}{\nu_\beta(B_{s-1})} \, < \, \frac12 \quad \mbox{and} \quad \sup_{s \ge 1} \frac{\nu_\beta(A_s \cap B_s)}{\nu_\beta(A_s)} \, < \, \frac12.$$ By symmetry of the argument, we just prove the second inequality above. Denote by $L$ the side length of the smallest cube $\mathrm{Q}_{2s}$ containing $A_s$. Then we have that $A\cap B$ can be decomposed into $\mathrm{C}_n = 8^n - 4^n$ cubes $S_j$ each of which satisfies  that $S_j = R_j + a_{S_j}$ for some cube $R_j = R_j(S_j) \subset A_s \setminus B_s$ of side length equal to $L/8$ and such that the angle between any point in $R_j$ and $a_{S_j}$ is smaller than $\pi/3$. We can also impose that $|a_{S_j}| \ge L/8$, see Figure II. This implies that for each $x$ in $R_j$ we have $$|x+a_{S_j}| \, \ge \, |x| + |a_{S_j}| \cos \sphericalangle (x, a_{S_j}) \, \ge \, |x| + \frac12 |a_{S_j}|.$$ Since $A_s \subset \R^n \setminus \mathrm{B}_1(0)$ for $s \ge 1$, we have
\begin{eqnarray*} 
\nu_{\beta}(A_s \cap B_s) & = & \int_{A_s \cap B_s} |x|^{-\beta} dx \ = \ \sum_{j=1}^{\mathrm{C}_n} \int_{S_j} |x|^{-\beta} dx \\ & = & \sum_{j=1}^{\mathrm{C}_n} \int_{R_j} |x+a_{S_j}|^{-\beta} dx \ \le \  \sum_{j=1}^{\mathrm{C}_n} \int_{R_j} \big( |x| + \frac12 |a_{S_j}| \big)^{-\beta} dx.
\end{eqnarray*}

\noindent
\begin{picture}(360,200)(-180,-100)
    \linethickness{1pt}
    \put(-22.5,-22.5){\line(0,1){45}}
    \put(22.5,-22.5){\line(0,1){45}}
    \put(-22.5,-22.5){\line(1,0){45}}
    \put(-22.5,22.5){\line(1,0){45}}
    
    \put(-45,-45){\line(0,1){90}}
    \put(45,-45){\line(0,1){90}}
    \put(-45,-45){\line(1,0){90}}
    \put(-45,45){\line(1,0){90}}
    
    \put(-90,-90){\line(0,1){180}}
    \put(90,-90){\line(0,1){180}}
    \put(-90,-90){\line(1,0){180}}
    \put(-90,90){\line(1,0){180}}    
    \linethickness{2pt}
    \put(-45.5,-46){\line(0,1){24.5}}
    \put(-22.5,-46){\line(0,1){24.5}}
    \put(-46,-45){\line(1,0){24.5}}
    \put(-46,-22.5){\line(1,0){24.5}}

    \put(-90,-46){\line(0,1){24.5}}
    \put(-67.5,-46){\line(0,1){24.5}}
    \put(-91,-45){\line(1,0){24.5}}
    \put(-91,-22.5){\line(1,0){24.5}}
    
    \put(22.5,21.5){\line(0,1){24.5}}
    \put(45,21.5){\line(0,1){24.5}}
    \put(21.5,22.5){\line(1,0){24.5}}
    \put(21.5,45){\line(1,0){24.5}}

    \put(22.5,44){\line(0,1){24.5}}
    \put(45,44){\line(0,1){24.5}}
    \put(21.5,45){\line(1,0){24.5}}
    \put(21.5,67.5){\line(1,0){24.5}} 
    
    \put(67.5,44){\line(0,1){24.5}}
    \put(90,44){\line(0,1){24.5}}
    \put(66.5,45){\line(1,0){24.5}}
    \put(66.5,67.5){\line(1,0){24.5}}         
     \linethickness{0.5pt}
 \qbezier[42](-90,67.5)(0,67.5)(90,67.5)
 \qbezier[42](-90,45)(0,45)(90,45)
 \qbezier[42](-90,22.5)(0,22.5)(90,22.5)
 \qbezier[42](-90,-22.5)(0,-22.5)(90,-22.5)
 \qbezier[42](-90,-45)(0,-45)(90,-45)
 \qbezier[42](-90,-67.5)(0,-67.5)(90,-67.5)

  \qbezier[42](67.5,-90)(67.5,0)(67.5,90)
  \qbezier[42](45,-90)(45,0)(45,90)
  \qbezier[42](22.5,-90)(22.5,0)(22.5,90)
  \qbezier[42](-22.5,-90)(-22.5,0)(-22.5,90)
  \qbezier[42](-45,-90)(-45,0)(-45,90)
  \qbezier[42](-67.5,-90)(-67.5,0)(-67.5,90)
  \qbezier[16](0,90)(0,60)(0,22.5)
  \qbezier[16](0,-90)(0,-60)(0,-22.5)
  \qbezier[16](90,0)(60,0)(22.5,0)
  \qbezier[16](-90,0)(-60,0)(-22.5,0)
    \linethickness{.2pt}
    \put(0,-100){\line(0,1){200}}
    \put(-100,0){\line(1,0){200}}
   \put(-84,-37){$S_1$}
   \put(-41,-37){$R_1$}
   \put(29,53){$S_2$}
   \put(74,53){$S_3$}
   \put(25,35){\tiny{$R_2$}}
   \put(28,28){\tiny{$= \! R_3$}}
   \put(4,-38){\large{$A_s \setminus B_s$}}
   \put(48,-82){\large{$A_s \cap B_s$}}

\end{picture}

\null

\vskip-25pt

\null

\begin{center}
\textsc{Figure II} \\ {\small There is a cube $R_j$ for each cube $S_j$}
\end{center}

\noindent Using that $|x| \leq \sqrt{n}L/2$ for $x \in A_s$  and $|a_{S_j}| \geq L/8$ $$\frac{1}{\nu_{\beta}(R_j)} \int_{R_j} \big( |x| + \frac12 |a_{S_j}| \big)^{-\beta} dx \, \le \, \sup_{x \in R_j} \frac{(|x| + \frac12 |a_{S_j}|)^{-\beta}}{|x|^{-\beta}} \, \le \, \left(\frac{\sqrt{n}}{\sqrt{n} + 1/8}\right)^{\beta}.$$ Therefore, we obtain
\begin{eqnarray*}
\frac{\nu_{\beta}(A_s \cap B_s)}{\nu_\beta(A_s)} & \le & \mathrm{C}_n \Big(\frac{\sqrt{n}}{\sqrt{n} + 1/8} \Big)^{\beta} \\ & \le & 8^n \Big(\frac{\sqrt{n}}{\sqrt{n} + 1/8} \Big)^{\beta} \ < \ \frac12 \quad \mbox{for} \quad \beta \gtrsim n^{\frac32}.
\end{eqnarray*}
A similar argument shows that $$\frac{\nu_\beta(A_s \cap B_{s-1})}{\nu_\beta(B_{s-1})} \, < \, \frac12$$ for $\beta \gtrsim n^{3/2}$ and $s \ge 2$, whereas the same estimate holds for $s=1$ as a consequence of the fact that $B_0$ contains $[-\lambda,\lambda]^n$ for $\lambda > 1$ large enough. This completes the construction of an admissible covering. It remains to construct filtrations $(\Sigma_{jk})_{k \ge 1}$ for $j=\mathrm{a}, \mathrm{b}$ which are regular and composed of doubling atoms. Recall that we set $\Sigma_{j1} = \Sigma_j$ and define $\Sigma_{j2}$ by splitting each atom in $\Sigma_j$ into a disjoint union of cubes. Namely, for $j=\mathrm{a}$ we keep $A_0$ and divide $A_s$ into the cubes $R_j, S_j$ in Figure II. We proceed similarly for $j=\mathrm{b}$. Once defined $\Sigma_{j2}$, we construct $\Sigma_{jk}$ by dyadic splitting of the cubes in $\Sigma_{j(k-1)}$. Note that the atoms in $\Sigma_{j1} \setminus \{A_0, B_0\}$ split at most into $8^n$ cubes $K$ centered at $c_K$ which are away from the origin. Thus $$\nu_{\beta}(2 K) \, = \, \int_{2K} |x|^{-\beta} dx \, \lesssim \, |2K| |c_K|^{-\beta} \, \lesssim \, |K| |c_K|^{-\beta} \, \lesssim \, \int_K |x|^{-\beta} dx \, = \, \nu_{\beta}(K). $$ It easily follows from this that all the atoms in $\Sigma_{jk}$ are doubling up to absolute constants independent of $k \ge 1$ and that both filtrations are regular. This shows that Theorem B1 applies to $(\R^n, \mu_\beta)$ with the Euclidean metric. 

\begin{remark}
\emph{A few comments are in order:}
\begin{itemize}
\item[i)] \emph{In the light of the example above, one could wonder what happens with the positive powers $d\mu_\gamma(x) = |x|^\gamma dx$ for $\gamma > 0$, but it is straightforward to show that these measures are doubling, so that we can handle them following the construction of Paragraph \ref{DCSubsection}.}

\vskip3pt

\item[ii)] \emph{Our proof of Theorem B1 relies crucially on the embedding of the space $\mathrm{DBMO}$ in $\mathrm{BMO}_{\Sigma_{\mathrm{ab}}}(\Omega)$ under suitable conditions. When the metric measure space $(\Omega,\mu,d)$ is of polynomial growth, we know from Tolsa \cite{To} that CZO's are $L_\infty \to \mathrm{RBMO}$ bounded. Since $\mathrm{RBMO} \subset \mathrm{DBMO}$, it is natural to wonder if we have $$\mathrm{RBMO} \subset \mathrm{BMO}_{\Sigma_{\mathrm{ab}}}(\Omega)$$ under weaker assumptions than in Theorem B1. It turns out that this is the case when there exists filtrations composed of doubling atoms, no matter whether they are regular or not. Indeed, noticing that RBMO can be described as a subspace of DBMO with an additional condition, it is this crucial extra condition introduced by Tolsa what allows to embed it into $\mathrm{BMO}_{\Sigma_{\mathrm{ab}}}(\Omega)$ and not in $\mathrm{bmo}_{\Sigma_{\mathrm{ab}}}(\Omega)$ for nonregular filtrations.} 
\end{itemize}
\end{remark}

\subsection{Concentration at the boundary}

Let $$d\mu_{\pm \alpha}(x) \, = \, e^{\pm |x|^\alpha} dx$$ on $\R^n$ equipped with the Euclidean metric. Carbonaro-Mauceri-Meda proved in \cite{CMM1,CMM2} that these measures satisfy their concentration condition when $\alpha > 1$. In this paragraph we shall prove that our hypothesis in Theorem B1 hold for any $\alpha > 0$, hence extending their results for measures with less concentration at the boundary. Let us start with the probability measure $\mu_{-\alpha}$. Pick $K = K(n,\alpha) > 0$ a large constant of the form $2^k$ for some $k \ge 1$ to be fixed below. Denote by $\mathbf{D}(\R^n)$ the standard filtration of dyadic cubes in $\R^n$. We consider the distinguished atom $A_0 = [-K,K]^n$. The other atoms $A_s \in \Pi_\mathrm{a}$ for $s \ge 1$ are chosen to be the cubes in $\mathbf{D}(\R^n \setminus A_0)$ which are maximal under the following constraint on the side-length $\ell(A_s)$ in terms of the modulus of its center $c_{A_s}$ $$\Pi_\mathrm{a} \, = \, \big\{ A_0 \big\} \bigcup \Big\{ A_s \mbox{ maximal in } \mathbf{D}(\R^n \setminus A_0) \ : \ \ell(A_s) \le K |c_{A_s}|^{1-\alpha}, \ s \ge 1 \Big\}.$$ Before defining $\Pi_\mathrm{b}$, we also need another dyadic filtration $\mathbf{D}'(\R^n)$ satisfying some specific properties which we now detail. Given cubes $(A,B) \in \mathbf{D}(\R^n) \times \mathbf{D}'(\R^n)$ of comparable size ---$2^{-\mathrm{k}_0} \le \ell(A)/\ell(B) \le 2^{\mathrm{k}_0}$ for some absolute constant $\mathrm{k}_0$--- with nonempty intersection, there exists a parallelepiped $R \subset A \triangle B$ such that: 

\begin{itemize}
\item[a)] $R$ is \lq\lq substantially closer" than $A \cap B$ to the origin,

\vskip2pt

\item[b)] There exists $a = a(R) \in \R^n$ such that $\displaystyle A \cap B \subset \bigcup_{j=1}^\mathrm{N} R + ja$,

\item[c)] $|a| \ge \frac{1}{\mathrm{N}} \max\{\ell(A), \ell(B)\}$ and $|x + ja| \ge |x| + \frac{1}{2} |a|$ for every $x \in R$.
\end{itemize}

\noindent
\begin{picture}(360,200)(-130,-100)
    \linethickness{.7pt}
    \put(30,-30){\line(0,1){75}}
    \put(105,-30){\line(0,1){75}}
    \put(30,-30){\line(1,0){75}}
    \put(30,45){\line(1,0){75}}
    
    \put(205,-130){\line(0,1){150}}
    \put(55,-130){\line(0,1){150}}
    \put(55,-130){\line(1,0){150}}
    \put(55,20){\line(1,0){150}}
    
     \qbezier[10](30,20)(45,20)(60,20)
     \qbezier[20](80,-30)(80,-5)(80,20)
    \linethickness{.7pt}
    \put(-100,-100){\line(0,1){180}}
    
    \put(-100,-100){\line(1,0){250}}
   \thicklines
   \qbezier(-100,-100)(-25,-56.25)(50,-12.5)
   \put(50,-12.5){\vector(1,0){25}}
   \linethickness{1pt}
   \put(47.5,-17.5){{\Huge $\cdot$}}
   \put(40,-24.5){$x$}
   \put(61,-19){$a$}
   \put(90,30){$A$}
   \put(155,-40){\huge $B$}
   \put(40,-5){\tiny{$R$}}
   \put(57,-5){\tiny{$R+a$}}
   \put(82,-5){\tiny{$R+2a$}}
\end{picture}

\vskip15pt 

\null

\begin{center}
\textsc{Figure III} \\ {\small $A \cap B$ is covered by at most $\mathrm{N}$ $a$-translates of $R$ \\ $R$ is \lq\lq substantially closer" to the origin than $A \cap B$ \\ $a$ is parallel to a coordinate axis for cubes $A, B$ in a sector around that axis \\ In particular, $x$ and $a$ are \lq\lq close" to being parallel, so that $|x + j a| \ge |x| + \frac12 |a|$}
\end{center}

\noindent Let $A_1$ be the cube in $\Pi_\mathrm{a} \setminus \{A_0\}$ whose center is the closest to the origin. Let $L = \ell(A_1)$ and pick $B_0 = A_0 + \frac13 L e_d$ with $e_d = (1,1, \ldots, 1)$. Then, the dyadic filtration $\mathbf{D}'(\R^n)$ is defined as one of the shifted dyadic filtrations in \cite{C} with initial cube being $B_0$. The fact that the properties above hold follows ultimately from the \lq\lq good separation" between $\mathbf{D}(\R^n)$ and $\mathbf{D}'(\R^n)$. Here we pick $K$ large enough so that the estimate $\mu_{-\alpha}(\frac12 A_0) > (1 - \varepsilon) \mu_{-\alpha}(\R^n)$ holds. In particular, we get $$\frac{\mu_{-\alpha}(A_0 \cap B_0)}{\mu_{-\alpha}(\R^n)} \, > \, 1 - \varepsilon$$ for some small $\varepsilon > 0$ to be fixed.  The family $\Pi_\mathrm{b}$ is defined similarly $$\Pi_\mathrm{b} \, = \, \big\{ B_0 \big\} \bigcup \Big\{ B_s \mbox{ maximal in } \mathbf{D}'(\R^n \setminus B_0) \, : \, \ell(B_s) \le K |c_{B_s}|^{1-\alpha}, \ s \ge 1 \Big\}.$$ Set $\Sigma_j = \sigma(\Pi_j)$ for $j=\mathrm{a}, \mathrm{b}$ and observe that $\Sigma_\mathrm{a} \cap \Sigma_\mathrm{b} = \{ \R^n, \emptyset \}$ by construction. Therefore, to prove that $(\Sigma_\mathrm{a}, \Sigma_\mathrm{b})$ yields an admissible covering we only need to check that we have $$\sup_{A \in \Pi_\mathrm{a} \setminus \{A_0\}} \sum_{B \in R_{A}} |R_B| \frac{\mu_{-\alpha}(A \cap B)^2}{\mu_{-\alpha}(A) \mu_{-\alpha}(B)} \, < \, 1.$$ According to our definition of $A_s$, it is a simple exercise to check that we have $\ell(A_s) \ge \frac13 K |c_{A_s}|^{1-\alpha}$ for all $s \ge 1$ but a finite number (independent of $K$) of cubes close to the origin. The same argument holds for atoms in $\Pi_\mathrm{b}$. In particular, we have $|R_A|, |R_B| \le \mathrm{C}_n$ for all $(A,B) \in \Pi_\mathrm{a} \times \Pi_\mathrm{b}$. Therefore, when $B = B_0$ we obtain $$|R_{B}| \frac{\mu_{-\alpha}(A\cap B)^2}{\mu_{-\alpha}(A)\mu_{-\alpha}(B)} \, \le \, \mathrm{C}_n \frac{\mu_{-\alpha}(A \cap B)}{\mu_{-\alpha}(B)} \, < \, \mathrm{C}_n \frac{\varepsilon}{1 - \varepsilon} \, < \, \frac12$$ for $\varepsilon < \frac13 \mathrm{C}_{n}^{-1}$. Otherwise, when $B \neq B_0$ we obtain  
\begin{eqnarray*}
\frac{\mu_{-\alpha}(A\cap B)}{\mu_{-\alpha}(R)} & = & \frac{1}{\mu_{-\alpha}(R)} \int_{A \cap B} e^{-|x|^{\alpha}} dx \\ & \le & \sum_{j=1}^\mathrm{N} \frac{1}{\mu_{-\alpha}(R)} \int_{R} e^{-|x+ja|^{\alpha}} dx \\ & \le & \frac{\mathrm{N}}{\mu_{-\alpha}(R)} \int_{R} e^{-(|x|+\frac12|a|)^{\alpha}} dx \ \le \ \mathrm{N} \sup_{x \in R} e^{-(|x| + \frac12 |a|)^{\alpha} + |x|^{\alpha}}. 
\end{eqnarray*}
If $\alpha =1$, we get an estimate $\mathrm{N} e^{- \frac12 |a|} \ge \mathrm{N} e^{- \mathrm{C}_n' K}$. For other values of $\alpha > 0$, an straightforward application of the mean value theorem gives $$\big( |x| + \frac12 |a| \big)^{\alpha} - |x|^{\alpha} \, \ge \, \frac{\alpha}{18} K |c_A|^{1-\alpha} |x|^{\alpha-1} \, \ge \, \mathrm{C}_n' K$$ since $|x| \sim |c_A|$. Hence we get for $A \neq A_0$
\begin{eqnarray*}
\sum_{B \in R_A} |R_B| \frac{\mu_{-\alpha}(A\cap B)^2}{\mu{-\alpha}(A)\mu{-\alpha}(B)} & < & \frac12 + \mathrm{C}_n^2 \sup_{B \neq B_0} \frac{\mu_{-\alpha}(A\cap B)^2}{\mu_{-\alpha}(A)\mu_{-\alpha}(B)} \\ & \le & \frac12 + \mathrm{C}_n^2 \mathrm{N} e^{- \mathrm{C}_n' K} \ < \ 1
\end{eqnarray*}
picking $K = K(n,\alpha)$ large enough. This shows the we have an admissible covering. Note that our choice of cubes for $\alpha > 1$ is a family which becomes smaller and smaller when we get away from the origin. This is in the spirit of Mauceri-Meda construction for the Gaussian measure \cite{MM}. On the contrary, when $\alpha < 1$ we pick larger and larger cubes as we get away from the origin. This construction seems not useful in \cite{CMM1,CMM2} since we may not use the locally doubling property for arbitrarily large cubes. Let us complete the proof by showing that the other hypotheses in Theorem B1 hold. Our choice of filtrations $(\Sigma_{jk})_{k \ge 1}$ for $j=\mathrm{a}, \mathrm{b}$ is by dyadic splitting of the cubes in $\Pi_\mathrm{a}$ and $\Pi_\mathrm{b}$ respectively. The regularity of such filtrations will follow from the fact that every atom in $(\Sigma_{jk})_{k \ge 1}$ is $(3,\beta)$-doubling for some absolute constant $\beta$, and this suffices to complete the proof. If $Q$ is any subcube of $A_0 \cup B_0$, there are dimensional constants $k_n$ and $K_n$ such that $$k_n |Q| \leq \int_Q e^{-|x|^{\alpha}} \; dx \leq K_n |Q|.$$
and hence $Q$ is trivially $(3,\beta)$-doubling. Otherwise, we compute
\begin{eqnarray*}
\frac{\mu_{-\alpha}(3Q)}{\mu_{-\alpha}(Q)} & \le & \frac{|3Q|}{|Q|} \sup_{x \in 3Q} e^{-|x|^{\alpha}} \sup_{x\in Q} e^{|x|^{\alpha}} \\
& \le & 3^n \exp \Big( \big( |x_Q| + \frac{1}{2} \sqrt{n} \ell(Q) \big)^{\alpha} - \big( |x_Q| - \frac{3}{2}\sqrt{n} \ell(Q) \big)^{\alpha} \Big) \ \le\ \beta
\end{eqnarray*}
for some absolute constant $\beta > 0$ using one more time the mean value theorem. 

\begin{remark}
\emph{A few comments are in order:}
\begin{itemize}
\item[i)] \emph{Given $\alpha > 0$ and by minor modifications in the above arguments, we may also produce an admissible covering for $(\R^n, e^{|x|^{\alpha}}dx)$ which satisfies the hypotheses of Theorem B1 with respect to the Euclidean metric.} 

\vskip3pt

\item[ii)] \emph{In this paper we $\wedge$-intersect two truncated martingale BMO spaces, but our results also hold for finite $\wedge$-intersections, details are simple and not very relevant. Mauceri-Meda BMO space for the Gaussian measure in \cite{MM} can be described as such a finite intersection of BMO spaces using a construction similar to the one above for $\mu_{-2}$ but intersecting $n+1$ BMO spaces instead of 2. Namely, one uses as many filtrations as needed to cover all cubes in $\R^n$ with dyadic cubes of comparable size, see for instance \cite{C} for the optimal choice. This establishes an inclusion of their BMO space into our 2-intersection $\mathrm{BMO}_{\Sigma_{\mathrm{ab}}}$ associated to $\mu_{-2}$, which still interpolates and it is strictly larger. The latter assertion can be proved following the argument which shows that classical BMO is strictly contained in dyadic BMO.}

\vskip3pt

\item[iii)] \emph{A geometric interpretation of our definition of admissible covering could be that we still impose certain concentration at the boundary, but much less than Carbonaro-Mauceri-Meda. In support of this, let us consider an admissible covering $(\Sigma_\mathrm{a}, \Sigma_\mathrm{b})$. Let $\A$ be a finite family in $\Pi_\mathrm{a} \setminus \{A_0\}$ and let $R_\A$ be the union $\cup_{A \in \A} R_A$. If we consider the set $R_\A$ as a measurable set and interpret $R_\A \setminus \A$ as the region \lq\lq close to the boundary\rq\rq${}$ then we can prove that $$\mu(R_\A) \, \le \, \frac{1}{1-c(\Sigma_{\mathrm{ab}})} \, \mu(R_\A \setminus \A)$$ or equivalently $\mu(\A) \le c(\Sigma_{\mathrm{ab}}) \, \mu(R_\A)$. Indeed, we have}
\begin{eqnarray*}
\hskip30pt \mu(\A) \!\!\! & = & \!\!\! \sum_{A \in \A} \sum_{B \in R_\A} \mu(A \cap B) \\ \!\!\! & \le & \!\!\! c(\Sigma_{\mathrm{ab}})^{\frac12} \sum_{A \in \A} \mu(A)^{\frac12} \Big( \sum_{B \in R_\A} \frac{\mu(B)}{|R_B|} \Big)^{\frac12} \\ \!\!\! & \le & \!\!\! c(\Sigma_{\mathrm{ab}})^{\frac12} \, \mu(\A)^{\frac12} \Big( \sum_{A \in \A} \sum_{B \in R_\A} \frac{\mu(B)}{|R_B|} \Big)^{\frac12} \, \le \, \big( c(\Sigma_{\mathrm{ab}}) \mu(\A) \mu(R_\A) \big)^{\frac12}.
\end{eqnarray*}
\end{itemize}
\end{remark}

\section{{\bf Calder\'on-Zygmund operators II}}

In this section we will study the class of atomic Calder\'on-Zygmund operators (ACZO) defined in the Introduction over a given measure space $(\Omega, \Sigma, \mu)$. More precisely, we shall prove Theorem B2 and illustrate it with a few constructions of dyadic operators satisfying its hypotheses. 

\demBB Following the same argument as in the proof of Theorem B1, we can use duality and our interpolation result in Theorem A to reduce the $L_p$-boundedness in the assertion to the $L_\infty(\Omega) \to \mathrm{BMO}_{\Sigma_{\mathrm{ab}}}(\Omega)$ boundedness of our  ACZO. This is however standard. Indeed, since the filtration is regular we know that $\mathrm{BMO}_{\Sigma_{\mathrm{ab}}}(\Omega) \simeq \mathrm{bmo}_{\Sigma_{\mathrm{ab}}}(\Omega)$. Up to absolute constants, the norm in the latter space is given by $$\|Tf\|_{\mathrm{bmo}_{\Sigma_{\mathrm{ab}}}(\Omega)} \, = \, \sup_{\mathrm{Q} \in \mathbf{\Pi}} \, \inf_{\mathrm{k_Q} \in \C} \Big( \frac{1}{\mu(\mathrm{Q})} \int_{\mathrm{Q}} \big| f(w) - \mathrm{k_Q} \big|^2 d\mu(\omega) \Big)^{\frac12},$$ where $\mathbf{\Pi} = \mathbf{\Pi}_\mathrm{a} \cup \mathbf{\Pi}_\mathrm{b}$ is the set of atoms in any of the two filtrations. Decompose $$f = f \chi_{\widehat{\mathrm{Q}}} + f \chi_{\Omega \setminus \widehat{\mathrm{Q}}} = f_1 + f_2.$$ As usual, we pick $\mathrm{k_Q} = (Tf_2)_{\mathrm{Q}}$. Then we control the term for $Tf_1$ using the $L_2$-boundedness of $T$ and the regularity of the filtrations. The term $Tf_2 - \mathrm{k_Q}$ is dominated by means of the H\"ormander kernel condition given in the definition of ACZO. Namely $$\Big( \frac{1}{\mu(\mathrm{Q})} \int_{\mathrm{Q}} \big| Tf_1(\omega) \big|^2 d\mu(\omega) \Big)^{\frac12} \, \le \, \|T\|_{2 \to 2} \, \sqrt{\frac{\mu(\widehat{Q})}{\mu(\mathrm{Q})}} \, \|f\|_\infty \, \lesssim \, \|f\|_\infty$$ by regularity of the filtrations. On the other hand 
\begin{eqnarray*}
\lefteqn{\Big( \frac{1}{\mu(\mathrm{Q})} \int_{\mathrm{Q}} \big| Tf_2(\omega) - (Tf_2)_{\mathrm{Q}} \big|^2 d\mu(\omega) \Big)^{\frac12}} \\ & \le & \Big( \frac{1}{\mu(\mathrm{Q})^2} \int_{\mathrm{Q}} \int_{\mathrm{Q}} \Big[ \int_{\Omega \setminus \widehat{\mathrm{Q}}} \big| K(z_1,x)- K(z_2,x) \big| d\mu(x) \Big]^2 d\mu(z_1) d\mu(z_2) \Big)^{\frac12} \|f\|_\infty
\end{eqnarray*}
which is dominated by $\|f\|_\infty$ according to the H\"ormander condition for ACZOs. \fin 

\begin{remark}
\emph{As mentioned in the Introduction, standard prototypes of atomic Calder\'on-Zygmund operator include martingale transforms, perfect dyadic CZOs and Haar shift operators. These are usually defined on the Euclidean space $\R^n$ equipped with a dyadic filtration. Nevertheless, the exact same arguments apply on any dyadically doubling measure space or even for any measure space equipped with a  two-sided regular filtration. In a recent paper \cite{LMP}, the authors studied those nondoubling measure spaces for which Haar shift operators satisfy weak type $(1,1)$ estimates. Theorem B2 provides a tool to produce nondoubling measure spaces over which Haar shifts or more general atomic CZOs are $L_\infty \to \mathrm{BMO}$ bounded. In the case of martingale transforms, Haar shift operators and perfect dyadic CZOs in $(\R^n,\mu)$, all of them satisfy the H\"ormander like condition $$\sup_{\mathrm{Q} \in \mathbf{D}(\R^n)} \, \sup_{z_1, z_2 \in \mathrm{Q}} \, \int_{\R^n \setminus \widehat{\mathrm{Q}}} \big| k(z_1,x) - k(z_2,x) \big| + \big| k(x,z_1) - k(x,z_2) \big| \, d\mu(x) \, < \, \infty.$$ This means that these operators are ACZOs satisfying Theorem B2 as long as we can find an admissible covering $(\Sigma_\mathrm{a}, \Sigma_\mathrm{b})$ and regular filtrations over it so that all the atoms are cubes in $\mathbf{D}(\R^n)$ or suitable unions of those. If we review our examples in Section \ref{Sect3}, this is not the case of our construction for $d\mu_{\pm \alpha}(x) = e^{\pm |x|^\alpha} dx$. It is however quite simple to adapt our construction for $$d\mu_\beta(x) = \frac{dx}{1 + |x|^\beta}$$ so that it satisfies the hypotheses of Theorem B2. In particular, the Haar shift operators defined on $(\R^n, \mu_\beta)$ are $L_\infty(\R^n) \to \mathrm{BMO}_{\Sigma_{\mathrm{ab}}}(\R^n)$ bounded. It remains open to decide whether an admissible covering exists on the exponential measure spaces $(\R^n, \mu_{\pm \alpha})$ using only atoms referred to one and not two dyadic systems.}
\end{remark}

\section{{\bf Matrix-valued forms of our results}}

In this section, we extend our main results to the context of operator-valued functions. Noncommutative forms of Calder\'on-Zygmund theory have been recently studied in \cite{JMP1,MP, Pa1}. There are however no specific results in the context of nondoubling metric measure spaces. Unfortunately, it seems difficult to extend the approach of Tolsa \cite{To} or Mauceri-Meda \cite{MM} to the operator-valued or even the noncommutative setting, since their interpolation results rest on good-$\lambda$ inequalities which do not have a noncommutative analogue so far. On the other hand, the semicommutative approach in \cite{Pa1} is valid for doubling spaces, but again presents serious obstacles to be extended to the nondoubling setting. The crucial aspect of our approach is that it ultimately rests on martingale inequalities which have been successfully transferred to the noncommutative setting. Namely, after Pixier-Xu's seminal contribution \cite{PX} on Burkholder-Gundy inequalities for noncommutative martingales, we find analogues of Doob's maximal inequalities; Gundy, Davis and atomic decompositions; Burkholder conditional square functions; John-Nirenberg inequalities; $L_p$/BMO interpolation results... see \cite{HM,J,JMe,JMu,JPe,JX,Mu,PR, Pe} and the references therein.     

Let us briefly introduce the framework for our results in this section, we refer to \cite[Section 1]{Pa1} for a rather complete review of the necessary background adapted to our necessities. We also refer the reader to Pisier-Xu's survey \cite{PX2} for more on noncommutative $L_p$ theory. Let $(\Omega,\Sigma,\mu)$ be a $\sigma$-finite measure space and consider any pair $(\M,\tau)$ given by a von Neumann algebra $\M$ equipped with a normal semifinite faithful trace $\tau$. This is sometimes called a noncommutative measure space. We will write $(\RR,\varphi)$ to denote the von Neumann algebra generated by essentially bounded functions $f: \Omega \to \M$ equipped with the trace $$\varphi(f) = \int_{\Omega} \tau(f(\omega)) \, d\mu(\omega).$$ $\RR$ is the von Neumann algebra tensor product $\RR = L_\infty(\Omega) \bar\otimes \M$ and we may consider the corresponding noncommutative spaces $L_p(\RR,\varphi)$. This semicommutative model is the context where we intend to generalize our main results. Apart from its own interest as an operator-valued model, it constitutes a first step towards further results for more general von Neumann algebras. In particular, as \cite{JMP0} demonstrates certain fully noncommutative questions can be reduced to the semicommutative setting. The reader which is not familiar with von Neumann algebra theory is encouraged to read this section restricting the attention to matrix-valued functions. In other words, replace $\M$ by the algebra $M_m$ of $m \times m$ matrices and $\tau$ by the standard trace $\mathrm{tr}$. The difficulties are similar in this case than in the general setting, as long as we provide results with constants independent of $m$. We also refer to \cite{Pa1} for a comparison between this model and the vector-valued setting, which differs substantially in the endpoint estimates. 

\subsection{The BMO spaces}

In this paragraph we review the definitions and results in Section \ref{Sect1} for the semicommutative setting described above. Given a filtration $(\Sigma_k)_{k \ge 1}$ of $(\Omega,\Sigma,\mu)$, we consider the conditional expectations $$f \in \RR \mapsto \mathsf{E}_{\Sigma_k} \otimes id_{\M}(f) \in \RR$$ still denoted by $\mathsf{E}_{\Sigma_k}$. The martingale $\mathrm{bmo}$ and $\mathrm{BMO}$ norms are 
\begin{eqnarray*}
\|f\|_{\mathrm{bmo}} & = & \max \Big\{ \|f\|_{\mathrm{bmo}_\mathrm{c}}, \|f^*\|_{\mathrm{bmo}_\mathrm{c}} \Big\}, \\
\|f\|_{\mathrm{BMO}} & = & \max \Big\{ \|f\|_{\mathrm{BMO}_\mathrm{c}}, \|f^*\|_{\mathrm{BMO}_\mathrm{c}} \Big\}, 
\end{eqnarray*}
where the column norms are defined as in the commutative case taking into account that we use $|x|^2 = x^*x$ for any operator $x$ on a Hilbert space. The interpolation result $[\mathrm{BMO}, L_1(\RR)]_{1/p} \simeq L_p(\RR)$ was proved by Musat in \cite{Mu} for any semifinite von Neumann algebra $\RR$. This is the noncommutative analogue of Janson-Jones's interpolation theorem. If we set $\|f\|_{\mathrm{h}_p^\mathrm{r}} = \|f^*\|_{\mathrm{h}_p^\mathrm{c}}$, where the norm in $\mathrm{h}_p^\mathrm{c}$ is defined as in the commutative case, then the noncommutative Hardy spaces have the form $$\mathrm{h}_p \, = \, \begin{cases} \mathrm{h}_p^\mathrm{r} + \mathrm{h}_p^\mathrm{c} & \mbox{if } 1\le p \le 2, \\ \mathrm{h}_p^\mathrm{r} \cap \mathrm{h}_p^\mathrm{c} & \mbox{if } 2 \le p \le \infty. \end{cases}$$ This combination of row and column square functions is known to be the right one for $L_p$ inequalities, as it was discovered for the first time with the noncommutative Khintchine inequalities \cite{LP,LPP}. The interpolation result $[\mathrm{bmo}, \mathrm{h}_1]_{1/p} \simeq \mathrm{h}_p$ was proved in \cite{BCP} for noncommutative martingales. As in the commutative case, the projections $J_{\Sigma_1} = id - \mathsf{E}_{\Sigma_1}$ are bounded on $\mathrm{bmo}$, $\mathrm{BMO}$, $L_p$ and $\mathrm{h}_p$, so that we will be working with these complemented subspaces which enjoy the same interpolation and duality properties than the original spaces. Note that the identity 
\begin{eqnarray*}
\|f\|_{J_{\Sigma_1}(\mathrm{bmo_c})} \!\!\! & = & \!\!\! \sup_{k \ge 1} \Big\| \Big( \mathsf{E}_{\Sigma_k} \big| f - \mathsf{E}_{\Sigma_k} f \big|^2 \Big)^\frac12 \Big\|_\M  \\ \!\!\! & = & \!\!\! \sup_{A \in \mathbf{\Pi}} \Big\| \Big( \frac{1}{\mu(A)} \int_A \Big| f(w) - \frac{1}{\mu(A)} \int_A f d\mu \Big|^2 d\mu(w) \Big)^\frac12 \Big\|_\M
\end{eqnarray*}
still holds and we have $J_{\Sigma_1}(\mathrm{bmo}) \simeq J_{\Sigma_1}(\mathrm{BMO})$ for regular filtrations. Consider an admissible covering $(\Sigma_{\mathrm{a}}, \Sigma_{\mathrm{b}})$ of $(\Omega,\Sigma,\mu)$ and any pair of martingale filtrations $(\Sigma_{jk})_{k \ge 1}$ with $\Sigma_{j1} = \Sigma_j$ for $j = \mathrm{a}, \mathrm{b}$. Denote by $\mathrm{BMO}_\mathrm{a}$ and $\mathrm{BMO}_\mathrm{b}$ the BMO spaces associated to this filtrations in the semicommutative algebra $\RR$ and set 
\begin{eqnarray*}
\mathrm{BMO}_{\Sigma_j}(\RR) & = & J_{\Sigma_j}(\mathrm{BMO}_j), \\ [5pt] \mathrm{BMO}_{\Sigma_\mathrm{ab}}(\RR) & = & \mathrm{BMO}_{\Sigma_\mathrm{a}}(\RR) \wedge \mathrm{BMO}_{\Sigma_\mathrm{b}}(\RR).
\end{eqnarray*}     
The John-Nirenberg inequalities, atomic descriptions of $\mathrm{H}_1$ and duality results have also been transferred to the context of noncommutative martingales \cite{BCP,HM,JMu,PX} and we will not review these results here, since they will not play a crucial role.  

\subsection{The interpolation theorem}

In this paragraph we state the analogue of Theorem A in the operator-valued setting. As usual we will write $L_p^\circ(\RR)$ for the subspace of mean-zero elements with respect to $\mu$. In the terminology we use for admissible coverings $$L_p^\circ(\RR) \, \simeq \, J_{\Sigma_\mathrm{a} \cap \Sigma_\mathrm{b}}(L_p(\RR)).$$    

\begin{theorem} \label{ThmASC}
Let $(\Sigma_\mathrm{a}, \Sigma_\mathrm{b})$ be an admissible covering in $(\Omega,\Sigma,\mu)$ and consider the semicommutative space $\RR = L_\infty(\Omega) \bar\otimes \M$. Then, for each $2 \le p < \infty$ there exists a constant $c_p \! \ge \! 1$ such that $$L_p^\circ(\RR) \, \simeq_{c_p} \, J_{\Sigma_\mathrm{a}}(L_p(\RR))  \, \wedge \, J_{\Sigma_\mathrm{b}}(L_p(\RR)).$$ In particular, we have by complex interpolation that $$\big[ \mathrm{BMO}_{\Sigma_{\mathrm{ab}}}(\RR), L_1^\circ (\RR) \big]_{1/q} \simeq_{c_q} L_q^\circ(\RR) \qquad (1 < q < \infty),$$ with $\mathrm{BMO}_{\Sigma_{\mathrm{ab}}}(\RR)$ constructed with any two martingale filtrations over  $(\Sigma_\mathrm{a}, \Sigma_\mathrm{b})$.
\end{theorem}

\sketch Thanks to the close connection with martingales, the proof is entirely parallel to the one given in the classical case. Indeed, combining standard facts from noncommutative $L_p$ theory with the martingale results reviewed in the previous paragraph, it is a simple exercise to adapt our proof of Theorem A to the present case. The only subtle point is the inequality $$\big\| \mathsf{E}_{\Sigma_\mathrm{a}} \mathsf{E}_{\Sigma_\mathrm{b}} f \big\|_p^p \, \le \, \big\| \mathsf{E}_{\Sigma_\mathrm{a}} \big| \mathsf{E}_{\Sigma_\mathrm{a}} f \big|^{\frac{p}{2}} \big\|_2^2,$$ which is used in the last two steps of our argument. Namely, in the classical case this is due to the conditional Jensen's inequality $\phi(\mathsf{E}_{\Sigma_k} f) \le \mathsf{E}_{\Sigma_k} \phi(f)$ for convex functions $\phi$. On the contrary, its noncommutative form does not hold for all $p \ge 2$ since we need the operator-convexity of the function $\phi(x) = |x|^\beta$ for $\beta = p/2$ and $x$ not necessarily positive. This is the case for $\beta \ge 2$ or equivalently $p \ge 4$, but it fails for $2 \le p < 4$. Note however that the ultimate goal in Steps 5 and 6 is to show that $\| \mathsf{E}_{\Sigma_\mathrm{a}} \mathsf{E}_{\Sigma_\mathrm{b}} f \|_p \le c_p(\Sigma_{\mathrm{ab}}) \|f\|_p$ for some $0 < c_p(\Sigma_{\mathrm{ab}}) < 1$. To prove it, we observe that 
\begin{eqnarray*}
\mathsf{E}_{\Sigma_\mathrm{a}}(g_1^* g_2) & = & \xi_k(g_1)^* \xi_k(g_2), \\ \mathsf{E}_{\Sigma_\mathrm{a}} \mathsf{E}_{\Sigma_\mathrm{b}} (f_1^* f_2) & = & \omega_k(f_1)^* \omega_k(f_2)
\end{eqnarray*}
for certain right $\RR_k$-module maps $\xi_k, \omega_k: L_q(\RR) \to C_q(L_q(\RR))$ with $\RR_k = \mathsf{E}_{\Sigma_k}(\RR)$. Namely, this follows from standard factorization properties of completely positive unital maps in terms of Hilbert modules, see for instance \cite{J}. Let  us consider the polar decompositions $f = u|f|$ and $g = v|g|$ of $g = \mathsf{E}_{\Sigma_\mathrm{b}}f$. Then we can factorize $\mathsf{E}_{\Sigma_\mathrm{a}} \mathsf{E}_{\Sigma_\mathrm{b}}f$ in two ways 
\begin{eqnarray*}
\mathsf{E}_{\Sigma_\mathrm{a}} g & = & \mathsf{E}_{\Sigma_\mathrm{a}} \big( v |g|^{\frac12} |g|^{\frac12} \big) \ = \ \xi_k \big( |g|^{\frac12}v^* \big)^* \xi_k \big( |g|^{\frac12} \big), \\ \mathsf{E}_{\Sigma_\mathrm{a}} \mathsf{E}_{\Sigma_\mathrm{b}}f & = & \mathsf{E}_{\Sigma_\mathrm{a}} \mathsf{E}_{\Sigma_\mathrm{b}} \big( u |f|^{\frac12} |f|^{\frac12} \big) \ = \ \omega_k \big( |f|^{\frac12}u^* \big)^* \omega_k \big( |f|^{\frac12} \big).
\end{eqnarray*}
This yields the following estimates
\begin{eqnarray*}
\hskip-35pt \big\| \mathsf{E}_{\Sigma_\mathrm{a}} \mathsf{E}_{\Sigma_\mathrm{b}} f \big\|_p & \le & \big\| \xi_k \big( |g|^{\frac12}v^* \big) \big\|_{2p} \big\| \xi_k \big( |g|^{\frac12} \big) \big\|_{2p} \\ & \le & \big\| \xi_k \big( |g|^{\frac12}v^* \big)^*\xi_k \big( |g|^{\frac12}v^* \big) \big\|_{p}^{\frac12} \big\| \xi_k \big( |g|^{\frac12} \big)^* \xi_k \big( |g|^{\frac12} \big) \big\|_{p}^{\frac12} \\ & = & \big\| \mathsf{E}_{\Sigma_\mathrm{a}} \big( v |g| v^* \big) \big\|_{p}^{\frac12} \big\| \mathsf{E}_{\Sigma_\mathrm{a}} \big( |g| \big) \big\|_{p}^{\frac12}  \, \le \, \big\| f \big\|_p^{\frac12} \big\| \mathsf{E}_{\Sigma_\mathrm{a}} \big| \mathsf{E}_{\Sigma_{\mathrm{b}}} f \big|^{\frac{p}{2}} \big\|_2^{\frac{1}{p}},
\end{eqnarray*}
and \begin{eqnarray*}
\big\| \mathsf{E}_{\Sigma_\mathrm{a}} \mathsf{E}_{\Sigma_\mathrm{b}} f \big\|_p & \le & \big\| \omega_k \big( |f|^{\frac12}u^* \big) \big\|_{2p} \big\| \omega_k \big( |f|^{\frac12} \big) \big\|_{2p} \\ & \le & \big\| \omega_k \big( |f|^{\frac12}u^* \big)^*\omega_k \big( |f|^{\frac12}u^* \big) \big\|_{p}^{\frac12} \big\| \omega_k \big( |f|^{\frac12} \big)^* \omega_k \big( |f|^{\frac12} \big) \big\|_{p}^{\frac12} \\ & = & \big\| \mathsf{E}_{\Sigma_\mathrm{a}} \mathsf{E}_{\Sigma_\mathrm{b}} \big( u |f| u^* \big) \big\|_{p}^{\frac12} \big\| \mathsf{E}_{\Sigma_\mathrm{a}} \mathsf{E}_{\Sigma_{\mathrm{b}}} \big( |f| \big) \big\|_{p}^{\frac12}  \, \le \, \big\| f \big\|_p^{\frac12} \big\| \mathsf{E}_{\Sigma_\mathrm{a}} \mathsf{E}_{\Sigma_{\mathrm{b}}} |f|^{\frac{p}{2}} \big\|_2^{\frac{1}{p}}.
\end{eqnarray*}
The last inequality in both estimates follows from Kadison-Schwarz  inequality for operator-convex functions, since $\phi(x) = x^\beta$ is operator-convex on $\R_+$ for $\beta \ge 1$. The first estimate is the right substitute in Step 5 and the second one in Step 6. \fin

\subsection{The Calder\'on-Zygmund operators}

We now consider Calder\'on-Zygmund operators in semicommutative algebras associated to operator-valued kernels. Our construction is standard, we refer to \cite{Du,JMP0,RRT} for further details. Let us we write $L_0(\M)$ for the $*$-algebra of $\tau$-measurable operators affiliated with $\M$ and consider kernels $k: (\Omega \times \Omega) \setminus \Delta \to \mathcal{L}(L_0(\M))$ defined away from the diagonal $\Delta $ of $\Omega \times \Omega$ and which take values in linear maps on $\tau$-measurable operators. If $d$ is a metric in $\Omega$, the standard H\"ormander kernel condition takes the same form in this setting when we replace the absolute value by the norm in the algebra $\mathcal{B}(\M)$ of bounded linear operators acting on $\M$ $$\sup_{\begin{subarray}{c} \mathrm{B} \, d\mathrm{-ball} \\ z_1, z_2 \in \mathrm{B} \end{subarray}} \, \int_{\Omega \setminus \alpha \mathrm{B}}  \big\| k(z_1,x) - k(z_2,x) \big\|_{\mathcal{B}(\M)} + \big\| k(x,z_1) - k(x,z_2) \big\|_{\mathcal{B}(\M)} \, d\mu(x) \, < \, \infty.$$ Define a CZO in $(\RR, \varphi,d)$ as any linear map $T$ satisfying the following properties 
\begin{itemize}
\item $T$ is bounded on $L_\infty(\M; L_2^r(\Omega))$ $$\hskip10pt \Big\| \Big( \int_\Omega Tf(x) Tf(x)^* \, d\mu(x) \Big)^{\frac12} \Big\|_\M \, \lesssim \, \Big\| \Big( \int_\Omega f(x) f(x)^* \, d\mu(x) \Big)^{\frac12} \Big\|_\M.$$

\item $T$ is bounded on $L_\infty(\M; L_2^c(\Omega))$ $$\hskip10pt \Big\| \Big( \int_\Omega Tf(x)^* Tf(x) \, d\mu(x) \Big)^{\frac12} \Big\|_\M \, \lesssim \, \Big\| \Big( \int_\Omega f(x)^* f(x) \, d\mu(x) \Big)^{\frac12} \Big\|_\M.$$

\item The kernel representation $$Tf(x) \, = \, \int_\Omega k(x,y) (f(y)) \, d\mu(y) \quad \mbox{holds for} \quad x \notin \mathrm{supp}_\Omega f$$ and some kernel $k: (\Omega \times \Omega) \setminus \Delta \to \C$ satisfying the H\"ormander condition.
\end{itemize}
The first two conditions replace the usual $L_2$-boundedness, see \cite{JMP0} for explanations.

\begin{theorem}
Let $(\Sigma_\mathrm{a},\Sigma_\mathrm{b})$ be an admissible covering of $(\Omega, \Sigma, \mu)$. Assume that $\Sigma$ admits regular filtrations $(\Sigma_{jk})_{k \ge 1}$ by successive refininement of $\Sigma_{j1} = \Sigma_j$ for $j=\mathrm{a},\mathrm{b}$ and that each atom in $\Sigma_{jk}$ is a $(\mathrm{C}_0, \alpha, \beta)$-doubling set for certain absolute constants $\mathrm{C}_0, \alpha, \beta > 0$. Let $\mathrm{BMO}_{\Sigma_{\mathrm{ab}}}(\RR)$ denote the $\wedge$-intersection of the $\mathrm{BMO}$ spaces defined over these filtrations. Then, every \emph{CZO} extends to a bounded map 
\begin{itemize}
\item[i)] $\mathrm{H}_{\Sigma_{\mathrm{ab}}}^1(\RR) \to L_1(\RR)$,

\item[ii)] $L_\infty(\RR) \to \mathrm{BMO}_{\Sigma_{\mathrm{ab}}}(\RR)$,

\item[iii)] $L_p^\circ(\RR) \to L_p^\circ(\RR)$ for $1 < p < \infty$.
\end{itemize}
Moreover, if $T$ is $L_2(\RR)$-bounded then $T: L_p(\RR) \to L_p(\RR)$ for all $1 < p < \infty$.
\end{theorem}

\dem According to Theorem \ref{ThmASC} (interpolation) and the semicommutative form of Remark \ref{H1Rem} (duality), it turns out that $L_\infty(\RR) \to \mathrm{BMO}_{\Sigma_{\mathrm{ab}}}(\RR)$ boundedness automatically implies $\mathrm{H}_{\Sigma_{\mathrm{ab}}}^1(\RR) \to L_1(\RR)$ boundedness, as well as $L_p^\circ(\RR) \to L_p^\circ(\RR)$ boundedness. Moreover, if $T$ is also $L_2$-bounded we may reproduce the argument given in the proof of Theorem B1 to obtain $L_p$-boundedness for all $1 < p < \infty$. Let us then focus from now on in the $L_\infty \to \mathrm{BMO}$ boundedness. Define $$\mathrm{DBMO} = \mathrm{DBMO}_\mathrm{r} \cap \mathrm{DBMO}_\mathrm{c}$$ with $\|f\|_{\mathrm{DBMO}_\mathrm{r}} = \|f^*\|_{\mathrm{DBMO_c}}$ and $$\|f\|_{\mathrm{DBMO_c}} \, = \, \sup_{\begin{subarray}{c} \mathrm{B} \, \mathrm{ball} \\ d\mathrm{-doubling} \end{subarray}} \Big\| \Big( \frac{1}{\mu(\mathrm{B})} \int_\mathrm{B} \Big| f(w) - \frac{1}{\mu(\mathrm{B})} \int_\mathrm{B} f \, d\mu \Big|^2 \, d\mu(w) \Big)^\frac12 \Big\|_\M.$$ As usual, we write $|x|^2$ for $x^*x$. The assertion follows from $$L_\infty(\RR) \stackrel{T}{\longrightarrow} \mathrm{DBMO} \stackrel{id}{\longrightarrow} \mathrm{bmo}_{\Sigma_{\mathrm{ab}}}(\RR) \simeq \mathrm{BMO}_{\Sigma_{\mathrm{ab}}}(\RR).$$ The boundedness of the chain above can be justified as in the proof of Theorem B1. Indeed, the analogies in the argument lead us to apply the new conditions which appear in our definition of semicommutative CZO, see \cite{JMP0}. \fin

\begin{remark}
\emph{Theorem B2 also admits an straightforward generalization to the semicommutative setting. Again, our use of martingale techniques makes the proof entirely analogous, so that we think it would be too repetitive to include it here.}
\end{remark}

\bibliographystyle{amsplain}

\begin{thebibliography}{99}
\bibitem {BL} J. Bergh and J. L\"{o}fstr\"{o}m, Interpolation Spaces. Springer-Verlag, Berlin, 1976.

\bibitem {BCP} T. Bekjan, Z. Chen, M. Perrin and Z. Yin, Atomic decomposition and interpolation for Hardy spaces of noncommutative martingales. J. Funct. Anal. \textbf{258} (2010), 2483-2505.

\bibitem{CMM1} A. Carbonaro, G. Mauceri, S. Meda, $H^1$ and BMO for certain locally doubling metric measure spaces. Ann. Scuola Normale Sup. Pisa \textbf{8} (2009), 543-582.

\bibitem{CMM2} A. Carbonaro, G. Mauceri, S. Meda, $H^1$ and BMO for certain locally doubling metric measure spaces of finite measure. Colloq. Math., \textbf{118} (2010), 13-41.

\bibitem {C} J.M. Conde, A note on dyadic coverings and nondoubling Calder\'on-Zygmund theory. J. Math. Anal. App. \textbf{397} (2013), 785-790.

\bibitem {CP} J.M. Conde-Alonso and J. Parcet, Atomic blocks for martingales. In progress. 

\bibitem {D} B. Davis, On the integrability of the martingale square function. Israel J. Math. \textbf{8} (1970), 187-190.

\bibitem {Du} J. Duoandikoetxea, Fourier Analysis. Translated and revised from the 1995 Spanish original by D. Cruz-Uribe. Grad. Stud. Math. \textbf{29}. American Mathematical Society, 2001.

\bibitem {G} A. Garsia, Martingale Inequalities, Seminar Notes on Recent Progress. Math. Lecture Notes Series 1973.

\bibitem {GJ} J.B. Garnett and P.W. Jones, BMO from dyadic BMO. Pacific J. Math. \textbf{99} (1982), 351-371.

\bibitem {HLMP} G. Hong, L.D. L\'opez-S\'anchez, J.M. Martell and J. Parcet, Calder\'on-Zygmund operators associated to matrix-valued kernels. Int. Math. Res. Not. \textbf{14} (2014), 1221-1252.

\bibitem {HM} G. Hong and T. Mei, John-Nirenberg inequality and atomic decomposition for noncommutative martingales. J. Funct. Anal. \textbf{263} (2011), 1064-1097. 

\bibitem {HK} T. Hyt\"onen and A. Kairema, Systems of dyadic cubes in a doubling metric space. Colloq. Math. \textbf{126} (2012), 1-33. 

\bibitem {JJ} S. Janson and P. Jones, Interpolation between $H_p$-spaces: the complex method. J. Funct. Anal. \textbf{48} (1982), 58-80.

\bibitem {JN} F. John and F. Nirenberg, On functions of bounded mean oscillation. Comm. Pure App. Math. \textbf{14} (1961), 415-426.

\bibitem {J} M. Junge, Doob's inequality for non-commutative martingales. J. Reine Angew. Math. \textbf{549} (2002), 149-190.

\bibitem{JMP0} M. Junge, T. Mei and J. Parcet, Smooth Fourier multipliers on group von Neumann algebras. ArXiv: 1010.5320. 

\bibitem{JMP1} M. Junge, T. Mei and J. Parcet, Algebraic Calder\'on-Zygmund theory. In progress. 

\bibitem {JMe} M. Junge and T. Mei, Noncommutative Riesz transforms -- A probabilistic approach. Amer. J. Math. \textbf{132} (2010), 611-681.

\bibitem {JMu} M. Junge and M. Musat, A noncommutative version of the John-Nirenberg theorem. Trans. Amer. Math. Soc. \textbf{359} (2007), 115-142.

\bibitem {JPe} M. Junge and M. Perrin, Theory of $\H_p$-spaces for continuous filtrations in von Neumann algebras. To appear in Asterisque.

\bibitem {JX} M. Junge and Q. Xu, Noncommutative Burkholder/Rosenthal inequalities. Ann. Probab. \textbf{31} (2003), 948-995.


\bibitem {LMP} L.D. L\'opez-S\'anchez, J.M. Martell and J. Parcet, Dyadic harmonic analysis beyond doubling measures. ArXiv: 1211.6291.

\bibitem {LP} F. Lust-Piquard, In\'{e}galit\'{e}s de Khintchine dans $C_p$ $(1 < p < \infty)$. C.R. Acad. Sci. Paris \textbf{303} (1986), 289-292.

\bibitem {LPP} F. Lust-Piquard and G. Pisier, Non-commutative Khintchine and Paley inequalities. Ark. Mat. \textbf{29} (1991), 241-260.

\bibitem {MMNO} J. Mateu, P. Mattila, A. Nicolau and J. Orobitg, BMO for nondoubling measures. Duke Math. J. \textbf{102} (2000), 533-565.

\bibitem{MM} G. Mauceri, S. Meda, BMO and $H^1$ for the Ornstein-Uhlenbeck operator. J. Funct. Anal. \textbf{252} (2007), 278-313.

\bibitem {M1} T. Mei, BMO is the intersection of two translates of dyadic BMO. C. R. Acad. Sci. Paris \textbf{336} (2003), 1003-1006.

\bibitem {M2} T. Mei, Operator Valued Hardy Spaces. Mem. Amer. Math. Soc. (2007) \textbf{881}.

\bibitem {MP} T. Mei and J. Parcet, Pseudo-localization of singular integrals and noncommutative Littlewood-Paley inequalities. Int. Math. Res. Not. \textbf{9} (2009), 1433-1487.

\bibitem {Mu} M. Musat, Interpolation between non-commutative BMO and non-commutative $L_p$-spaces. J. Funct. Anal. \textbf{202} (2003), 195-225. 

\bibitem {NTV1} F. Nazarov, S. Treil and A. Volberg, Accretive system $Tb$-theorems on nonhomogeneous spaces. Duke Math. J. \textbf{113} (2002), 259-312.

\bibitem {Pa1} J. Parcet, Pseudo-localization of singular integrals and noncommutative Calder\'on-Zygmund theory. J. Funct. Anal. \textbf{256} (2009), 509-593.

\bibitem {PR} J. Parcet and N. Randrianantoanina, Gundy's decomposition  for non-commutative martingales and applications. Proc. London Math. Soc. \textbf{93} (2006), 227-252.

\bibitem {Pe} M. Perrin, A noncommutative Davis' decomposition for martingales. J. London Math. Soc. \textbf{80} (2009), 627-648.



\bibitem {PX} G. Pisier and Q. Xu, Non-commutative martingale inequalities. Comm. Math. Phys. \textbf{189} (1997), 667-698.

\bibitem {PX2} G. Pisier and Q. Xu, Non-commutative $L_p$-spaces. Handbook of the Geometry of Banach Spaces II (Eds. W.B. Johnson and J. Lindenstrauss) North-Holland (2003), 1459-1517.

\bibitem {RRT} J.L. Rubio de Francia, F. Ruiz and J.L. Torrea, Calder\'on-Zygmund theory for operator-valued kernels. Adv. Math. \textbf{62} (1986), 7-48.

\bibitem {To} X. Tolsa, BMO, $H^1$, and Calder\'on-Zygmund operators for non doubling measures. Math. Ann. \textbf{319} (2001), 89-149.
\end{thebibliography}


\vskip40pt




\hfill \noindent \textbf{Jose M. Conde-Alonso} \\
\null \hfill Instituto de Ciencias Matem{\'a}ticas \\ \null \hfill
CSIC-UAM-UC3M-UCM \\ \null \hfill Consejo Superior de
Investigaciones Cient{\'\i}ficas \\ \null \hfill C/ Nicol\'as Cabrera 13-15.
28049, Madrid. Spain \\ \null \hfill\texttt{jose.conde@icmat.es}

\vskip5pt

\hfill \noindent \textbf{Tao Mei} \\
\null \hfill Department of Mathematics
\\ \null \hfill Wayne State University \\
\null \hfill 656 W. Kirby Detroit, MI 48202. USA \\
\null \hfill\texttt{mei@wayne.edu}

\vskip5pt

\enlargethispage{2cm}

\hfill \noindent \textbf{Javier Parcet} \\
\null \hfill Instituto de Ciencias Matem{\'a}ticas \\ \null \hfill
CSIC-UAM-UC3M-UCM \\ \null \hfill Consejo Superior de
Investigaciones Cient{\'\i}ficas \\ \null \hfill C/ Nicol\'as Cabrera 13-15.
28049, Madrid. Spain \\ \null \hfill\texttt{javier.parcet@icmat.es}
\end{document}